\documentclass[final]{siamltex}

\usepackage{amssymb} \usepackage{cite} \usepackage{graphicx}
\usepackage{psfrag} \usepackage{subfigure} \usepackage{url}
\usepackage{amsmath} \usepackage{color}

\usepackage[colorlinks=false,urlcolor=blue,citecolor=blue,linkcolor=blue,
  bookmarks=true, bookmarksopen=false,pdftitle=created by
  dvipdf,pdfcreator=NM,
  pdfauthor=Megiddo,pdfsubject=mor.sty:6.29.04]{hyperref}

\def\a{\alpha}

\newcommand{\remove}[1]{} 
\newcommand{\dyn}[1]{} 

\newcommand{\EXP}[1]{\mathsf{E}\!\left[#1\right] }
\newcommand{\Var}[1]{\mathsf{Var}\!\left[#1\right] }

\newcommand{\rem}[1]{} 
\newtheorem{assumption}{Assumption}

\def\l{\left} \def\r{\right}
\def\nn{\nonumber}  \def\liminf{\mathop{\rm liminf}}
\def\limsup{\mathop{\rm limsup}}

\title{Distributed Stochastic Subgradient Projection Algorithms for
        Convex Optimization}

\author{S.~Sundhar~Ram, A.~Nedi\'{c}, and V.~V.~Veeravalli \thanks{The
        first and the third authors are with the Electrical and
        Computer Engineering Department at the University of Illinois at
        Urbana-Champaign.  The second author is with the Industrial
        and Enterprise Systems Engineering Department at 
        University of Illinois at
        Urbana-Champaign. They can be contacted at {\texttt
        \{ssriniv5,angelia,vvv\}@illinois.edu}. This work has been supported 
        by NSF Career Grant CMMI 07-42538.}}

\begin{document}
\maketitle

\begin{abstract}
We consider a distributed multi-agent network system where the goal is
to minimize a sum of convex objective functions of the agents subject
to a common convex constraint set.  Each agent maintains an iterate
sequence and communicates the iterates to its neighbors.  Then, each
agent combines weighted averages of the received iterates with its own
iterate, and adjusts the iterate by using subgradient information
(known with stochastic errors) of its own function and by projecting
onto the constraint set.

The goal of this paper is to explore the effects of stochastic
subgradient errors on the convergence of the algorithm.  We first
consider the behavior of the algorithm in mean, and then the
convergence with probability 1 and in mean square.  We consider
general stochastic errors that have uniformly bounded second moments
and obtain bounds on the limiting performance of the algorithm in mean
for diminishing and non-diminishing stepsizes.  When the means of the
errors diminish, we prove that there is mean consensus between the
agents and mean convergence to the optimum function value for
diminishing stepsizes.  When the mean errors diminish sufficiently
fast, we strengthen the results to consensus and convergence of the
iterates to an optimal solution with probability 1 and in mean square.

\end{abstract}



\pagestyle{myheadings} \thispagestyle{plain} \markboth{\ \ \
  S.~Sundhar~Ram, A.~Nedi\'{c} and V.~V.~Veeravalli}{Distributed
  Stochastic Subgradient Algorithms for Convex Optimization\ \ \ }

\section{Introduction}
\label{sec:intro}
A number of problems that arise in the context of wired and wireless
networks can be posed as the minimization of a sum of functions, when
each component function is available only to a specific agent
\cite{Rabbat05,Sundhar08a,Sundhar08d}. Often, it is not efficient, or
not possible, for the network agents to share their objective
functions with each other or with a central coordinator. In such
scenarios, distributed algorithms that only require the agents to
locally exchange limited and high level information are preferable.
For example, in a large wireless network, energy is a scarce resource
and it might not be efficient for a central coordinator to learn the
individual objective functions from each and every agent
\cite{Rabbat05}.  In a network of databases from which information is
to be mined, privacy considerations may not allow the sharing of the
objective functions \cite{Park02}.  In a distributed network on a
single chip, for the chip to be fault tolerant, it is desirable to
perform the processing in a distributed manner to account for the
statistical process variations\cite{Varatkar07}.

We consider constrained minimization of a sum of convex functions,
where each component function is known partially (with stochastic
errors) to a specific network agent.  The algorithm proposed builds on
the distributed algorithm proposed in \cite{Nedic08a} for the
unconstrained minimization problem.  Each agent maintains an iterate
sequence and communicates the iterates to 
its neighbors.  Then, each agent averages the received
iterates with its own iterate, and  adjusts the iterate by using
subgradient information (known with stochastic errors) of its own
function and by projecting onto the constraint set.  The inter-agent
information exchange model is a synchronous and delayless version of
the computational model proposed by Tsitsiklis \cite{Tsitsiklis84}.
The algorithm is \emph{distributed} since there is no central
coordinator. The algorithm is {\it local} since each agent uses only
locally available information (its objective function) and
communicates locally with its immediate neighbors.

Related to this work are the distributed incremental algorithms, where
the network agents 
sequentially update an iterate sequence in a
cyclic or a random order
\cite{Blatt07,Rabbat05,Nedic01,Sundhar08a,Johansson07}.  The effects of
stochastic errors on these algorithms have been investigated in
\cite{Bertsekas00,Gaivoronski94,Kiwiel03,Solodov98b,Nedic07,
Rabbat05,Sundhar08a}.  In an incremental algorithm, there is a single
iterate sequence and only one agent updates the iterate at a given
time. Thus, while being distributed and local, incremental algorithms differ
fundamentally from the algorithm studied in this paper (where all
agents update simultaneously).  Also related are the optimization
algorithms in \cite{Bertsekas97,Tsitsiklis86}. However, these
algorithms are not local as the complete objective function
information is available to each and every agent, with the
aim of distributing the processing.

The work in this paper is also related at a much broader level to the
distributed consensus algorithms
\cite{Tsitsiklis84,Tsitsiklis86,Bertsekas97,
Kar07,Olfati05,Xiao07,Lobel08,Nedic08b,Nedic08c,Olshevsky08,Jadbabaie03}. In
these algorithms, each agent starts with a different value and through
local information exchange, the agents eventually 
agree on a common value.  
The effect of random errors on consensus algorithms
have been investigated in \cite{Huang07,Lobel08,Xiao07,Kar07}.  In addition,
since we are interested in the effect of stochastic errors, our paper
is also related to the literature on stochastic subgradient methods
\cite{Ermoliev76,Ermoliev83,Ermoliev88}.

We consider general stochastic errors that have uniformly bounded
second moments and obtain bounds on the limiting performance of the
algorithm in mean for diminishing and non-diminishing stepsizes. 
When the means of the errors diminish, we prove
that there is mean consensus between the
agents and mean convergence to the optimum function value 
for diminishing stepsizes.  When 
the mean errors diminish sufficiently fast, we strengthen
the results to consensus and convergence of the iterates 
to an optimal solution with probability 1 and in mean square.

Our work expands the multi-agent distributed optimization framework
studied in \cite{Nedic08a}. The new contributions are: 1) the study of the
effects of stochastic errors in subgradient evaluations; 2) the
consideration of constrained optimization problem within the
distributed multi-agent setting.  The presence of the constraint
set complicates the analysis as it introduces non-linearities in the
system dynamics. The non-linearity issues that we face have some
similarities to those in the constrained consensus problem
investigated in \cite{Nedic08c}, though the problems are fundamentally
different. The presence of subgradient stochastic errors adds another
layer of complexity to the analysis as the errors made by each agent
propagate through the network to every other agent and also across
time, making the iterates statistically dependent across time and
agents.

The rest of the paper is organized as follows.  In
Section~\ref{sec:assum}, we formulate the problem, describe the
algorithm and state our basic assumptions. In Section~\ref{sec:prelim}
we state some results from literature that we use in the
analysis, while in Section~\ref{sec:bir}, we derive two important lemmas
that form the backbone of the analysis. 
In Section~\ref{sec:generror}, we study the convergence properties 
of the method in mean, and in Section~\ref{sec:dimerror}
we focus on  the convergence properties with probability 1 and in mean square.
Finally, we discuss some implications 
and provide some
concluding remarks in Sections~\ref{sec:implications} and 
\ref{sec:conclusion}.

\section{Problem, algorithm and assumptions}
\label{sec:assum}
In this section, we formulate the problem of interest and describe the
algorithm that we propose. We also state and discuss our assumptions
on the agent connectivity and information exchange. 

\subsection{Problem}
We consider a network of $m$ agents that are indexed by $1,\ldots, m$.
Often, when convenient, we index the agents by using set $V=
\{1,\ldots,m\}$. The network objective is to solve the following
constrained optimization problem:
\begin{align}
  \mbox{minimize \ \ } & \sum_{i=1}^{m} f_i(x) \cr 
  \mbox{subject to \ } & x \in X, \label{eqn:problem}
\end{align}
where $X\subseteq\Re^n$ is a constraint set and $f_i:X \to \Re$
for all $i$. Related to the problem %(\ref{eqn:problem}), 
we use the following notation
\[
f(x) = \sum_{i=1}^{m} f_i(x), \qquad f^* = \min_{x \in X} f(x), \qquad
X^* = \{x \in X: f(x) = f^*\}.
\]

We are interested in the case when the problem in (\ref{eqn:problem}) is 
convex. Specifically, we assume that the following assumption holds.
\begin{assumption}
  The functions $f_i$ and the set $X$ are such that
  \label{ass:convex}
  \begin{enumerate}
  \item[(a)] The set $X$ is closed and convex.  
  \item[(b)] The functions $f_i,\ i\in V$ are defined and convex 
over an open set that contains the set $X$. 
  \end{enumerate}
  \end{assumption}

The function $f_i$ is known only partially to agent $i$ in the sense
that the agent can only obtain a noisy estimate of the function
subgradient.  The goal is to solve problem (\ref{eqn:problem}) using
an algorithm that is distributed and local.\footnote{See
\cite{Sundhar08a,Sundhar08d} for wireless network applications that
can be cast in this framework.}

We make no assumption on the differentiability of the
functions $f_i.$ At points where the gradient does not exist, we use
the notion of subgradients. A vector $\nabla f_i$ is a subgradient of
$f_i$ at a point $x\in{\rm dom\, f}$ if the following relation holds
\begin{eqnarray}
\nabla f_i(x)^T (y - x)\leq f_i(y) - f_i(x) \qquad\hbox{for all $y\in
{\rm dom}\, f$}.
\label{eqn:lip3}
\end{eqnarray}
Since the set $X$ is contained in an open set over which the functions
are defined and convex, 
a subgradient of $f_i$ exists at any point of the set $X$
(see \cite{Bertsekas03} or \cite{Rockafellar70}).  

\subsection{Algorithm}

To solve the problem in (\ref{eqn:problem}) with its inherent 
decentralized information access,
we consider an iterative subgradient method. The iterations
are distributed accordingly among the agents, whereby
each agent $i$ is minimizing its convex objective 
$f_i$ over the set $X$ and locally exchanging the iterates 
with its neighbors.

Let $w_{i,k}$ be the iterate with agent~$i$ at the end of iteration
$k.$ At the beginning of iteration $k+1$, agent $i$ receives the
current iterate of a subset of the agents. Then, agent $i$ computes a
weighted average of these iterates and adjusts this average along the
negative subgradient direction of $f_i$, which is computed with
stochastic errors.  The adjusted iterate is then projected onto the
constraint set $X.$ Mathematically, each agent $i$ generates its
iterate sequence $\{w_{i,k}\}$ according to the following relation:
\begin{align}
  w_{i,k+1} = P_{X} \l[ v_{i,k} - \alpha_{k+1} \left(\nabla
    f_i\l(v_{i,k}\r) + \epsilon_{i,k+1} \right) \r],
  \label{eqn:sagrad}
\end{align}
starting with some initial iterate $w_{i,0}\in X.$ Here, $\nabla
f_i\l(v_{i,k}\r)$ denotes the subgradient of $f_i$ at $v_{i,k}$ and
$\epsilon_{i,k+1}$ is the stochastic error in the subgradient
evaluation. The scalar $\alpha_{k+1}>0$ is the stepsize and $P_{X}$
denotes the Euclidean projection onto the set $X.$ The vector
$v_{i,k}$ is the weighted average computed by agent $i$ and is given
by
\begin{equation}
v_{i,k} = \sum_{j\in N_i(k+1)} a_{i,j}(k+1) w_{j,k},
\label{eqn:averi}
\end{equation}
where $N_{i}(k+1)$ denotes the set of agents whose current iterates
are available to agent $i$ in the $(k+1)$-st iteration. We assume that
$i \in N_i(k+1)$ for all agents and at all times $k$.  The scalars
$a_{i,j}(k+1)$ are the non-negative weights that agent $i$ assigns to
agent $j$'s iterate. We will find it convenient to define
$a_{i,j}(k+1)$ as $0$ for $j \notin N_i(k+1)$ and rewrite
(\ref{eqn:averi}) as
\begin{equation}
v_{i,k} = \sum_{j=1}^{m} a_{i,j}(k+1) w_{j,k}.
\label{eqn:averig}
\end{equation}
This is  a ``consensus''-based step ensuring that, in a long run,
the information 
of each $f_i$ reaches every agent with the same frequency,
directly or through a sequence of 
local communications.
Due to this, the iterates $w_{j,k}$ become eventually ``the same'' 
for all $j$ and for large enough $k$.
The update step in (\ref{eqn:sagrad}) is just a subgradient iteration
for minimizing $f_i$ over $X$ taken after the ``consensus''-based step.

\subsection{Additional assumptions}
In addition to Assumption~\ref{ass:convex}, we make some
assumptions on the inter-agent exchange model and the weights. The
first assumption requires the agents to communicate sufficiently often
so that all the component functions, directly or indirectly, influence
the iterate sequence of any agent. Recall that we defined $N_i(k+1)$
as the set of agents that agent $i$ communicates with in iteration
$k+1.$ Define $(V,E_{k+1})$ to be the graph with edges
\[E_{k+1} = \{(j,i): j \in N_{i}(k+1), \, i\in V\}.\]
\begin{assumption}
\label{ass:con}
 There exists a scalar $Q$ such that the graph $(V,
\cup_{l=1,\ldots,Q} E_{k+l} )$ is strongly connected for all $k.$ 
\end{assumption}

It is also important that the influence of the functions $f_i$ is
``equal'' in a long run so that the sum of the component functions is
minimized rather than a weighted sum of them.  The influence of a
component $f_j$ on the iterates of agent $i$ depends on the weights
that agent $i$ uses.  To ensure equal influence, we make the following
assumption on the weights.
\begin{assumption}\label{ass:weight}
 For  $i \in V$ and all $k,$
  \begin{enumerate}
  \item [(a)] $a_{i,j}(k+1)\ge0$, and $a_{i,j}(k+1) = 0 $ when $j \notin
        N_{i}(k+1),$
  \item [(b)]$\sum_{j=1}^{m} a_{i,j}(k+1) = 1,$
  \item [(c)]There exists a scalar $\eta,$ $0 < \eta < 1,$ such that
    $a_{i,j}(k+1) \geq \eta$ when $j \in N_{i}(k+1),$
  \item [(d)]$\sum_{i=1}^{m} a_{i,j}(k+1) = 1.$
  \end{enumerate}
\end{assumption}

Assumptions~\ref{ass:weight}a and \ref{ass:weight}b state that each
agent calculates a weighted average of all the iterates it has access
to. Assumption~\ref{ass:weight}c ensures that each agent gives a
sufficient weight to its current iterate and all the iterates it
receives.\footnote{The agents need not be aware of the common bound
$\eta.$} Assumption~\ref{ass:weight}d, together with
Assumption~\ref{ass:con}, as we will see later, ensures that all the
agents are equally influential in the long run. 
In other words,  Assumption~\ref{ass:weight}d is crucial 
to ensure that $\sum_{i=1}^m f_i$ is minimized as opposed to a weighted 
sum of the functions $f_i$ with non-equal weights. To satisfy
Assumption~\ref{ass:weight}d, the agents need to coordinate their
weights. Some coordination schemes are discussed in
\cite{Nedic08a,Sundhar08a}.

\section{Preliminaries}
\label{sec:prelim}
In this section, we state some results for future reference.
\subsection{Euclidean norm inequalities}
For any vectors $v_1,\ldots,v_M \in \Re^n,$ we have
\begin{align}
  \sum_{i=1}^{M} \l\| v_i - \frac{1}{M} \sum_{j=1}^{M} v_j \r\|^2 \leq
  \sum_{i=1}^{M} \l\| v_i - x \r\|^2 \qquad \mbox{ for any $x \in
  \Re^n$}.
  \label{eqn:ineq3} 
\end{align}
The preceding relation states that the average of a finite set of
vectors minimizes the sum of distances between each vector and any
vector in $\Re^n$, which can be verified using the first-order 
optimality conditions. 

Both the Euclidean norm and its square are convex functions, i.e.,
for any vectors $v_1, \ldots, v_M \in\Re^n$ and
nonnegative scalars $\beta_1, \ldots,\beta_M$ such that
$\sum_{i=1}^{M} \beta_i = 1,$ we have 
\begin{align}
  \l\|\sum_{i=1}^{M} \beta_i v_i \r\| &\leq \sum_{i=1}^{M} \beta_i \|
  v_i \|, \label{eqn:ineq5} \\ \l\|\sum_{i=1}^{M} \beta_i v_i \r\|^2
  &\leq \sum_{i=1}^{M} \beta_i \| v_i \|^2. \label{eqn:ineq6}
\end{align}
The following inequality is the well-known\footnote
{See  for example \cite{Bertsekas03}, Proposition 2.2.1.}
{\it non-expansive} property of the
Euclidean projection onto a
nonempty, closed and convex set $X$, 
\begin{align}
  \|P_{X}[x] - P_{X}[y] \| \leq \| x - y \|
  \qquad \mbox{ for all $x,y \in \Re^n$}.
  \label{eqn:proj1}
\end{align}
%This property is known as the {\it non-expansive} property of the
%Euclidean projection (see \cite{Bertsekas03}, Proposition 2.2.1).

\subsection{Scalar sequences}
For a scalar $\beta$ and a scalar sequence $\{\gamma_k\}$, we consider
the ``convolution'' sequence
$\sum_{\ell=0}^k\beta^{k-\ell}\gamma_\ell =\beta^k\gamma_0
+\beta^{k-1}\gamma_1+\cdots+\beta\gamma_{k-1}+\gamma_k.$ We have the
following result.
\begin{lemma}
  \label{lemma:betalpha}
Let $\{\gamma_k\}$ be a scalar sequence.
\begin{itemize} 
\item [(a)] If $\lim_{k\to\infty}\gamma_k=\gamma$  and $0< \beta<1,$ then
$\lim_{k\to\infty}\sum_{\ell=0}^k \beta^{k-\ell}\gamma_\ell
=\frac{\gamma}{1-\beta}.$
\item [(b)] If $\gamma_k\ge0$ for all $k$,
 $\sum_k\gamma_k<\infty$ and $0< \beta<1,$  then
  $\sum_{k=0}^\infty
  \l(\sum_{\ell=0}^k\beta^{k-\ell}\gamma_\ell\r)<~\infty.$
\item [(c)] 
If $\limsup_{k\to\infty}\gamma_k=\gamma$  
and $\{\zeta_k\}$ is a positive scalar sequence
with $\sum_{k=1}^{\infty} \zeta_k =~\infty,$ then 
$\limsup_{K\to\infty}\frac{\sum_{k=0}^{K} \gamma_k \zeta_k
  }{\sum_{k=0}^{K} \zeta_k} \le \gamma.$
In addition, if $\liminf_{k\to\infty}\gamma_k=\gamma$, then 
$\lim_{K\to\infty}\frac{\sum_{k=0}^{K} \gamma_k \zeta_k
  }{\sum_{k=0}^{K} \zeta_k}~=~\gamma.$
\end{itemize}
\end{lemma}
 
\begin{proof} (a)\ 
Let $\epsilon>0$ be arbitrary.  Since $\gamma_k\to\gamma$ and 
for all $k$, there is an
index $K$ such that $|\gamma_k-\gamma|\le\epsilon$ for all $k\ge K$. For all
$k\ge K+1$, we have
\[\sum_{\ell=0}^k\beta^{k-\ell}\gamma_\ell=
\sum_{\ell=0}^K\beta^{k-\ell}\gamma_\ell +
\sum_{\ell=K+1}^k\beta^{k-\ell}\gamma_\ell \le \max_{0\le t\le
K}\gamma_t \sum_{\ell=0}^K\beta^{k-\ell}
+(\gamma+\epsilon) \sum_{\ell=K+1}^k\beta^{k-\ell}.\] 
Since
$\sum_{\ell=K+1}^k\beta^{k-\ell}\le\frac{1}{1-\beta}$ and
\[\sum_{\ell=0}^K\beta^{k-\ell}=\beta^k+\cdots+\beta^{k-K}
=\beta^{k-K}(1+\cdots+\beta^K)\le \frac{\beta^{k-K}}{1-\beta},\]
it follows that for all $k\ge K+1$,
\[\sum_{\ell=0}^k\beta^{k-\ell}\gamma_\ell
\le \left(\max_{0\le t\le K}\gamma_t\right)\frac{\beta^{k-K}}{1-\beta}
+\frac{\gamma+ \epsilon}{1-\beta}.\]
Therefore,
\[\limsup_{k\to\infty}\sum_{\ell=0}^k\beta^{k-\ell}\gamma_\ell
\le \frac{\gamma+\epsilon}{1-\beta}.\]
Since $\epsilon$ is arbitrary, we conclude that
$\limsup_{k\to\infty}\sum_{\ell=0}^k\beta^{k-\ell}\gamma_\ell
\le \frac{\gamma}{1-\beta}$.

Similarly, we have 
\[\sum_{\ell=0}^k\beta^{k-\ell}\gamma_\ell
\ge \min_{0\le t\le
K}\gamma_t \sum_{\ell=0}^K\beta^{k-\ell}
+(\gamma-\epsilon) \sum_{\ell=K+1}^k\beta^{k-\ell}.\] 
Thus, 
\[\liminf_{k\to\infty}
\sum_{\ell=0}^k\beta^{k-\ell}\gamma_\ell
\ge \liminf_{k\to\infty} \left(\min_{0\le t\le
K}\gamma_t \sum_{\ell=0}^K\beta^{k-\ell}
+(\gamma-\epsilon) \sum_{\ell=K+1}^k\beta^{k-\ell}\right).\]
Since 
$\sum_{\ell=0}^K\beta^{k-\ell}\ge \beta^{k-K}$
%=\beta^k+\cdots+\beta^{k-K}=\beta^{k-K}(1+\cdots+\beta^K)
and $\sum_{\ell=K+1}^k\beta^{k-\ell}=
\sum_{s=0}^{k-(K+1)}\beta^s$, which tends to $1/(1-\beta)$ as $k\to\infty$, 
it follows that
\[\liminf_{k\to\infty}\sum_{\ell=0}^k\beta^{k-\ell}\gamma_\ell
\ge \left(\min_{0\le t\le K}\gamma_t\right) \lim_{k\to\infty}\beta^{k-K}
+(\gamma-\epsilon)\lim_{k\to\infty}\sum_{s=0}^{k-(K+1)}\beta^s
=\frac{\gamma-\epsilon}{1-\beta}.\] 
Since $\epsilon$ is arbitrary, we have
$\liminf_{k\to\infty}\sum_{\ell=0}^k\beta^{k-\ell}\gamma_\ell
\ge \frac{\gamma}{1-\beta}.$
This and the relation 
$\limsup_{k\to\infty}\sum_{\ell=0}^k\beta^{k-\ell}\gamma_\ell\le 
\frac{\gamma}{1-\beta}$,
imply 
\[\lim_{k\to\infty}\sum_{\ell=0}^k\beta^{k-\ell}\gamma_\ell
=\frac{ \gamma}{1-\beta}.\]

\vskip 0.1 pc
\noindent
(b)\ Let $\sum_{k=0}^\infty \gamma_k<\infty.$
For any integer $M\ge1$, we have
\[\sum_{k=0}^M \l(\sum_{\ell=0}^k\beta^{k-\ell}\gamma_\ell\r)
=\sum_{\ell=0}^M\gamma_\ell\sum_{t=0}^{M-\ell}\beta^t
\le \sum_{\ell=0}^M\gamma_\ell\frac{1}{1-\beta},\]
implying that 
\[\sum_{k=0}^\infty \l(\sum_{\ell=0}^k\beta^{k-\ell}\gamma_\ell\r)
\le \frac{1}{1-\beta}\sum_{\ell=0}^\infty\gamma_\ell<\infty.\]

\vskip 0.1 pc
\noindent
(c)\ Since $\limsup_{k\to\infty} 
\gamma_k =\gamma,$ for every $\epsilon > 0$ there
is a large enough $K$ such that $\gamma_k \le \gamma+\epsilon$ 
for all $k \geq K.$ Thus, for any $M>K,$
\begin{align*}
  \frac{\sum_{k=0}^{M} \gamma_k \zeta_k }{\sum_{k=0}^{M}
    \zeta_k} &= \frac{\sum_{k=0}^{K} \gamma_k \zeta_k
    }{\sum_{k=0}^{M} \zeta_k} + \frac{\sum_{k=K+1}^{M}
    \gamma_k \zeta_k }{\sum_{k=0}^{M} \zeta_k}  
    \leq
    \frac{\sum_{k=0}^{K} \gamma_k \zeta_k }{\sum_{k=0}^{M}
    \zeta_k} + (\gamma+\epsilon) \frac{\sum_{k=K+1}^{M} \zeta_k
    }{\sum_{k=0}^{M} \zeta_k}. 
\end{align*}
By letting $M\to\infty$ and using $\sum_{k} \zeta_k = \infty,$ 
we see that 
$\limsup_{M\to\infty}\ \frac{\sum_{k=0}^{M} \gamma_k \zeta_k }{\sum_{k=0}^{M}
    \zeta_k} \le\gamma+\epsilon,$
and since $\epsilon$ is arbitrary, the result for the limit superior follows.

Analogously, if $\liminf_{k\to\infty} 
\gamma_k =\gamma,$ then for every $\epsilon > 0$ there
is a large enough $K$ such that $\gamma_k \ge \gamma-\epsilon$ 
for all $k \geq K.$ Thus, for any $M>K,$
\begin{align*}
  \frac{\sum_{k=0}^{M} \gamma_k \zeta_k }{\sum_{k=0}^{M}
    \zeta_k} &= \frac{\sum_{k=0}^{K} \gamma_k \zeta_k
    }{\sum_{k=0}^{M} \zeta_k} + \frac{\sum_{k=K+1}^{M}
    \gamma_k \zeta_k }{\sum_{k=0}^{M} \zeta_k}  
    \ge
    \frac{\sum_{k=0}^{K} \gamma_k \zeta_k }{\sum_{k=0}^{M}
    \zeta_k} + (\gamma-\epsilon) \frac{\sum_{k=K+1}^{M} \zeta_k
    }{\sum_{k=0}^{M} \zeta_k}. 
\end{align*}
Letting $M\to\infty$ and using $\sum_{k} \zeta_k = \infty,$ 
we obtain 
$\liminf_{M\to\infty}\ \frac{\sum_{k=0}^{M} \gamma_k \zeta_k }{\sum_{k=0}^{M}
    \zeta_k} \ge\gamma-\epsilon.$ Since $\epsilon>0$ is arbitrary, we have
$\liminf_{M\to\infty}\ \frac{\sum_{k=0}^{M} \gamma_k \zeta_k }{\sum_{k=0}^{M}
    \zeta_k} \ge\gamma$. This relation and the relation 
for the limit superior yield 
$\lim_{M\to\infty}\ \frac{\sum_{k=0}^{M} \gamma_k \zeta_k }{\sum_{k=0}^{M}
    \zeta_k} =\gamma$ when $\gamma_k\to\gamma$.
\end{proof}

\subsection{Matrix convergence}
Let $A(k)$ be the matrix with $(i,j)$-th entry equal to $a_{i,j}(k).$
As a consequence of Assumptions~\ref{ass:weight}a,~\ref{ass:weight}b
and~\ref{ass:weight}d, the matrix $A(k)$ is doubly
stochastic\footnote{The sum of its entries in every row and in every
column is equal to $1$.}.  Define, for all $k,s$ with $k \geq s,$
\begin{align}
\Phi(k,s) = A(k) A(k-1) \cdots A(s+1). \label{eqn:phi}
\end{align}
We next state a result from \cite{Nedic08b} (Corollary 1) on the
convergence properties of the matrix $\Phi(k,s).$ Let
$[\Phi(k,s)]_{i,j}$ denote the $(i,j)$-th entry of the matrix
$\Phi(k,s),$ and let $e\in\Re^m$ be the column vector with all entries
equal to $1.$
\begin{lemma}
  \label{lemma:lit}
  Let Assumptions~\ref{ass:con} and \ref{ass:weight} hold. Then
  \begin{enumerate}
    \item $\lim_{k \to \infty} \Phi(k,s) = \frac{1}{m}\, e e^T$ for all
      $s.$
    \item Further, the convergence is geometric and the rate of
    convergence is given by
      \begin{align*}
	\l| [\Phi(k,s)]_{i,j} - \frac{1}{m} \r| \leq \theta
	\beta^{k-s},
	%\label{eqn:roc}
      \end{align*}
      where
      \[
      \theta = \l(1 - \frac{\eta}{4 m^2} \r)^{-2} \qquad
      \beta = \l(1 - \frac{\eta}{4 m^2} \r)^{\frac{1}{Q}}.
      \]
  \end{enumerate}
\end{lemma}

\subsection{Stochastic convergence}
We next state some results that deal with the convergence of a
sequence of random vectors.  The first result is the well known
Fatou's lemma \cite{Billingsley79}.
\begin{lemma}
\label{lemma:Fat}
Let $\{X_i\}$ be a sequence of non-negative random variables. Then
$$ \EXP{\liminf_{n \to \infty} X_n} \leq \liminf_{n \to \infty}
\EXP{X_n}.
$$
\end{lemma}

The next
result is due to Robbins and Siegmund (Lemma~11,
Chapter~2.2, \cite{Polyak87}).
\begin{theorem}
\label{lemma:asm}
Let $\{B_k\}, \{D_k\},$ and $\{H_k\}$ be non-negative random sequences
and let $\{\zeta_k\}$ be a deterministic nonnegative scalar sequence. 
Let $G_k$ be the
$\sigma-$algebra generated by $B_1, \ldots, B_k, D_1, \ldots, D_k,H_1,
\ldots, H_k.$ Suppose that $\sum_{k} \zeta_k < \infty,$
\begin{equation}
\EXP{B_{k+1} \mid G_k } \leq (1 + \zeta_k) B_{k} - D_{k} + H_{k}\qquad
\hbox{for all $k$,} \label{eqn:asm}
\end{equation}
and $\sum_{k} H_k <\infty$ with probability 1.  Then, the sequence
$\{B_{k}\}$ converges to a non-negative random variable and $\sum_{k}
D_{k} < \infty$ with probability~1, and in mean.
\end{theorem}

\section{Basic relations}
\label{sec:bir}
In this section, we derive two basic relations that form the basis for
the analysis in this paper. 
The first of them deals with the disagreements
among the agents, and the second deals with the agent iterate sequences.

\subsection{Disagreement Estimate}
The agent disagreements are typically thought of as
the norms $\|w_{i,k}-w_{j,k}\|$ of  the differences between the iterates
$w_{i,k}$ and $w_{j,k}$   
generated by different agents according to 
(\ref{eqn:sagrad})--(\ref{eqn:averi}). 
Alternatively, the agent disagreements can be measured with 
respect to a reference sequence, which we adopt here.
In particular, we study the
behavior of $\|y_k-w_{i,k}\|,$ where $\{y_k\}$ is the auxiliary vector
sequence defined by
\begin{align}
  y_{k} &= \frac{1}{m} \sum_{i=1}^{m} w_{i,k}\qquad\hbox{for all $k$}. 
\label{eqn:ydef}
\end{align}
In the next lemma, we provide a basic estimate for $\|y_k-w_{j,k}\|.$
The rate of convergence result from Lemma~\ref{lemma:lit}
plays a crucial role in obtaining this estimate. 

\begin{lemma}
  \label{lemma:yandwj}
  Let Assumptions~\ref{ass:convex}a, \ref{ass:con}, and
  \ref{ass:weight} hold. Assume that the subgradients of $f_i$ are
  uniformly bounded over the set $X$, i.e., there are scalars $C_i$
  such that
   \[\|\nabla f_i(x)\|\le C_i\qquad\hbox{for all $x\in X$ and all $i\in V$}.
   \]  
   Then, for all $j\in V$ and $k\ge0$,
    \begin{align*}
      \|y_{k+1} - w_{j,k+1}\| \leq& m \theta \beta^{k+1} \max_{i\in V}
      \l\|w_{i,0}\r\| + \theta \sum_{\ell=1}^{k} \alpha_{\ell}
      \beta^{k+1-\ell} \sum_{i=1}^{m} \l( C_i + \|\epsilon_{i,\ell}\|
      \r) \\ &+ \frac{\alpha_{k+1}}{m} \sum_{i=1}^{m} (C_i +
      \|\epsilon_{i,k+1}\|) + \alpha_{k+1} (C_j +
      \|\epsilon_{j,k+1}\|).
    \end{align*} 
\end{lemma}

\begin{proof}
  Define for all $i\in V$ and all $k,$
  \begin{align}
    p_{i,k+1} &= w_{i,k+1} - \sum_{j=1}^{m} a_{i,j}(k+1)
    w_{j,k}. \label{eqn:pdef}
  \end{align}
  Using the matrices $\Phi(k,s)$ defined in (\ref{eqn:phi}) we can
  write %(as in \cite{Nedic07})
  \begin{align}
    w_{j,k+1} = \sum_{i=1}^{m} [\Phi(k+1,0)]_{j,i} w_{i,0} + p_{j,k+1}
    + \sum_{\ell=1}^{k}\l(\sum_{i=1}^{m} [\Phi(k+1,\ell)]_{j,i} p_{i,\ell}
    \r). \label{eqn:aform2}
  \end{align}
  Using (\ref{eqn:pdef}), we can also rewrite $y_{k},$ defined in
  (\ref{eqn:ydef}), as follows
  \begin{align*}
    y_{k+1} &= \frac{1}{m} \l( \sum_{i=1}^{m} \sum_{j=1}^{m}
    a_{i,j}(k+1) w_{j,k} + \sum_{i=1}^{m} p_{i,k+1} \r) \\ 
    &=\frac{1}{m} \l( \sum_{j=1}^{m} \l(\sum_{i=1}^{m} a_{i,j}(k+1) \r)
    w_{j,k} + \sum_{i=1}^{m} p_{i,k+1} \r). 
    \end{align*}
  In the view of the doubly stochasticity of the weights, we have
  $\sum_{i=1}^{m} a_{i,j}(k+1) =1$, implying that
  \begin{align*}
    y_{k+1} = \frac{1}{m} \l(
    \sum_{j=1}^{m} w_{j,k} + \sum_{i=1}^{m} p_{i,k+1} \r) = y_{k}
    + \frac{1}{m} \sum_{i=1}^{m} p_{i,k+1}.
  \end{align*}
  Therefore
  \begin{align}
    y_{k+1} &= y_{0} + \frac{1}{m} \sum_{\ell=1}^{k+1}
    \sum_{i=1}^{m} p_{i,\ell}  
    = \frac{1}{m}\sum_{i=1}^{m}
    w_{i,0} + \frac{1}{m} \sum_{\ell=1}^{k+1} \sum_{i=1}^{m}
    p_{i,\ell}. \label{eqn:astopf2}
  \end{align}

  Substituting for $y_{k+1}$ from (\ref{eqn:astopf2}) and for
  $w_{j,k+1}$ from (\ref{eqn:aform2}), we obtain {\small
     \begin{eqnarray*}
      \|y_{k+1} - w_{j,k+1}\| &=& \l\| \frac{1}{m}\sum_{i=1}^{m}
      w_{i,0} + \frac{1}{m} \sum_{\ell=1}^{k+1} \sum_{i=1}^{m}
      p_{i,\ell} \r. \nn \cr
      && \l. - \l( \sum_{i=1}^{m}
      [\Phi(k+1,0)]_{j,i} w_{i,0} + p_{j,k+1} +
      \sum_{\ell=1}^{k}\sum_{i=1}^{m} [\Phi(k+1,\ell)]_{j,i}
      p_{i,\ell} \r) \r\| \nn \\
      &=& \l\| \sum_{i=1}^{m} \l(\frac{1}{m}
      - [\Phi(k+1,0)]_{j,i} \r) w_{i,0} \r. \nn \\
      && \l. +
      \sum_{\ell=1}^{k}\sum_{i=1}^{m} \l(\frac{1}{m} -
      [\Phi(k+1,\ell)]_{j,i} \r) p_{i,\ell} + \l(
      \frac{1}{m}\sum_{i=1}^{m} p_{i,k+1} - p_{j,k+1} \r) \r\|.
  \end{eqnarray*}}
  Therefore, for all $j\in V$ and all $k$, 
    {\small
    \begin{align*}
      \|y_{k+1} - w_{j,k+1}\| \leq& \sum_{i=1}^{m} \l|\frac{1}{m} -
      [\Phi(k+1,0)]_{j,i} \r| \l\|w_{i,0} \r\| \\&+
      \sum_{\ell=1}^{k} \sum_{i=1}^{m} \l|\frac{1}{m} -
      [\Phi(k+1,\ell)]_{j,i} \r| \l\|p_{i,\ell} \r\| +
      \frac{1}{m}\sum_{i=1}^{m} \l\| p_{i,k+1} \r\| + \l\|p_{j,k+1}
      \r\|.
  \end{align*} }

  We can bound $\l\|w_{i,0} \r\| \leq \max_{i\in V}\|w_{i,0}\|.$
  Further, we can use the rate of convergence result from
  Lemma~\ref{lemma:lit} to bound $\l|\frac{1}{m} -
  [\Phi(k,\ell)]_{j,i} \r|.$ We obtain
   \begin{align}
      \|y_{k+1} - w_{j,k+1}\| \leq& m \theta \beta^{k+1}
      \max_{i\in V}\|w_{i,0}\| + \theta \sum_{\ell=1}^{k}
      \beta^{k+1-\ell} \sum_{i=1}^{m} \l\|p_{i,\ell} \r\| \cr
      &+
      \frac{1}{m}\sum_{i=1}^{m} \l\| p_{i,k+1} \r\| +
      \l\|p_{j,k+1} \r\|.
      \label{eqn:ywjp}
  \end{align} 

  We next estimate the norms of the vectors $\|p_{i,k}\|$ for any $k$.
  From the definition of $p_{i,k+1}$ in (\ref{eqn:pdef}) and the
  definition of the vector $v_{i,k}$ in~(\ref{eqn:averi}), we have
  $p_{i,k+1}=w_{i,k+1}-v_{i,k}$. Note that, being a convex combination
  of vectors $w_{j,k}$ in the convex set $X$, the vector $v_{i,k}$ is
  in the set $X$. By the definition of the iterate $w_{i,k+1}$
  in~(\ref{eqn:sagrad}) and the non-expansive property of the
  Euclidean projection in~(\ref{eqn:proj1}), we have
     \begin{align*} 
       \| p_{i,k+1} \| &= \l\|P_{X} \l[ v_{i,k} -
        \alpha_{k+1}\l(\nabla f_i(v_{i,k}) + \epsilon_{i,k+1}\r) \r] -
        v_{i,k}\r\| \\ &\leq \alpha_{k+1} \l\| \nabla f_i(v_{i,k}) +
        \epsilon_{i,k+1}\r\| \\ &\leq \alpha_{k+1} \l(C_i +
        \|\epsilon_{i,k+1}\| \r).
     \end{align*}
     In the last step we have used the subgradient boundedness. By
     substituting the preceding relation in (\ref{eqn:ywjp}), we
     obtain the desired relation.
\end{proof}

\subsection{Iterate Relation}
Here, we derive a 
relation for the distances $\|v_{i,k+1}-~z\|$ and the function value 
differences $f(y_k)-f(z)$ for an arbitrary $z\in X.$ 
This relation together with Lemma~\ref{lemma:yandwj}
provides the basis for our subsequent convergence analysis.
In what follows, recall that $f=\sum_{i=1}^m f_i$.

\begin{lemma}
   \label{lemma:biralpha}
   Let Assumptions~\ref{ass:convex}, \ref{ass:con}, and \ref{ass:weight}
   hold.  Assume that the subgradients of $f_i$ are uniformly bounded
   over the set $X$, i.e., there are scalars $C_i$ such that
   \[\|\nabla f_i(x)\|\le C_i\qquad\hbox{for all $x\in X$ and all $i\in V$}.\] 
   Then, for any $z \in X$ and all $k$,
   \begin{eqnarray*}
     \sum_{i=1}^{m} \|v_{i,k+1} - z\|^2 &\leq& \sum_{i=1}^{m} \|
     v_{i,k} - z\|^2 -2 \alpha_{k+1}\, \l( f(y_{k}) - f(z) \r) \cr
     &&+2\alpha_{k+1}\l(\max_{i\in V} C_i\r)\sum_{j=1}^m
     \|y_k-w_{j,k}\|\cr &&-2\alpha_{k+1} \sum_{i=1}^{m}
     \epsilon_{i,k+1}^T(v_{i,k}-z) + \alpha_{k+1}^2 \sum_{i=1}^{m} \l(
     C_i + \|\epsilon_{i,k+1}\| \r)^2.
   \end{eqnarray*}
 \end{lemma}
\begin{proof}
  Using the Euclidean projection property in (\ref{eqn:proj1}), from
  the definition of the iterate $w_{i,k+1}$ in (\ref{eqn:sagrad}), we
  have for any $z\in X$ and all $k$,
  \begin{align*}
    \|w_{i,k+1} - z\|^2 =& \l\| P_{X}\l[v_{i,k} - \alpha_{k+1} \l(
      \nabla f_i(v_{i,k}) + \epsilon_{i,k+1} \r) \r] - z \r\|^2 \\
      \leq& \|v_{i,k} - z\|^2 - 2\alpha_{k+1} \nabla f_i(v_{i,k})^T
      (v_{i,k} - z) - 2\alpha_{k+1} \epsilon_{i,k+1}^T(v_{i,k} - z)
      \\& + \alpha_{k+1}^2 \l\| \nabla f_i(v_{i,k}) + \epsilon_{i,k+1}
      \r\|^2.
  \end{align*}
  By using the subgradient inequality in (\ref{eqn:lip3}) to bound the
  second term, we obtain
  \begin{align}
    \|w_{i,k+1} - z\|^2 \leq& \|v_{i,k} - z\|^2 - 2\alpha_{k+1} \l(
    f_i(v_{i,k}) - f_i(z) \r) \cr&- 2\alpha_{k+1}
    \epsilon_{i,k+1}^T(v_{i,k} - z) + \alpha_{k+1}^2 \l\| \nabla
    f_i(v_{i,k}) + \epsilon_{i,k+1} \r\|^2.
  \label{eqn:sumvi}
  \end{align}
  
  Note that by the convexity of the squared norm
  [cf.~Eq.~(\ref{eqn:ineq6})], we have
  \[\sum_{i=1}^{m} \|v_{i,k+1} - z\|^2 = \sum_{i=1}^{m}
  \l\|\sum_{j=1}^{m} a_{i,j}(k+2) w_{j,k+1} - z \r\|^2 \leq
  \sum_{i=1}^{m} \sum_{j=1}^{m} a_{i,j}(k+2) \| w_{j,k+1} - z\|^2. \]
  In view of Assumption~\ref{ass:weight}, we have $\sum_{i=1}^{m}
  a_{i,j}(k+2)=1$ for all $j$ and $k$, implying that
  \[\sum_{i=1}^{m} \|v_{i,k+1} - z\|^2 
  \leq
  \sum_{j=1}^{m} \| w_{j,k+1} - z\|^2.\]
  
  By summing the relations in (\ref{eqn:sumvi}) over all $i\in V$ and
  by using the preceding relation, we obtain
  \begin{eqnarray}
    \sum_{i=1}^{m} \|v_{i,k+1} - z\|^2 &\leq& \sum_{i=1}^{m} \|
    v_{i,k} - z\|^2 - 2 \alpha_{k+1} \sum_{i=1}^{m} \l( f_i(v_{i,k}) -
    f_i(z) \r) \cr &&- 2 \alpha_{k+1} \sum_{i=1}^{m}
    \epsilon_{i,k+1}^T(v_{i,k} - z) + \alpha_{k+1}^2 \sum_{i=1}^{m}
    \l\|\nabla f_i(v_{i,k}) + \epsilon_{i,k+1}\r\|^2.
  \qquad\label{eqn:sumviz}
  \end{eqnarray}
  From (\ref{eqn:lip3}) we have
  \begin{align}
    f_i(v_{i,k}) - f_i(z) \geq& \l( f_i(v_{i,k}) - f_i(y_k) \r) +
    \l( f_i(y_{k}) - f_i(z) \r)\cr\ \geq& - \|\nabla f_i(v_{i,k})\|
    \|y_k - v_{i,k}\| + \l( f_i(y_{k}) - f_i(z)
    \r). \label{eqn:rcon}
  \end{align}
  Recall that $v_{i,k}=\sum_{j=1}^m a_{i,j}(k+1)w_{j,k}$
  [cf.~(\ref{eqn:averig})].  Substituting for $v_{i,k}$ and using the
  convexity of the norm [cf.~(\ref{eqn:ineq5})], from (\ref{eqn:rcon})
  we obtain
   \begin{align*}
     \sum_{i=1}^{m} f_i(v_{i,k}) - f_i(z) \geq& - \sum_{i=1}^{m}
      \|\nabla f_i(v_{i,k})\| \|y_k - v_{i,k}\| + \l( f(y_{k}) - f(z)
      \r) \nn \\ \geq& - \sum_{i=1}^{m} \|\nabla f_i(v_{i,k})\|
      \l\|y_k - \sum_{j=1}^{m} a_{i,j}(k+1) w_{j,k}\r\| + \l( f(y_{k})
      - f(z) \r)\nn \\ \geq& - \sum_{i=1}^{m} \|\nabla f_i(v_{i,k})\|
      \sum_{j=1}^{m} a_{i,j}(k+1) \|y_k - w_{j,k}\| + \l( f(y_{k}) -
      f(z) \r) \nn \\ \geq& - \l(\max_{i \in V} \|\nabla
      f_i(v_{i,k})\| \r) \sum_{j=1}^{m} \l(\sum_{i=1}^{m} a_{i,j}(k+1)
      \r) \|y_k - w_{j,k}\| \nn\\ &+ \l( f(y_{k}) - f(z) \r)\nn \\ =&
      - \l(\max_{i \in V} \|\nabla f_i(v_{i,k})\| \r) \sum_{j=1}^{m}
      \|y_k - w_{j,k}\| + \l( f(y_{k}) - f(z) \r).
      %\label{eqn:rcon2}
   \end{align*}
   By using the preceding estimate in relation (\ref{eqn:sumviz}),
   %and   (\ref{eqn:rcon2}), 
   we have
   \begin{align*}
     \sum_{i=1}^{m} \|v_{i,k+1} - z\|^2 \leq& \sum_{i=1}^{m} \|
     v_{i,k} - z\|^2 - 2\alpha_{k+1} \l( f(y_{k})- f(z) \r) \cr &+ 2
     \alpha_{k+1} \l(\max_{i \in V} \|\nabla f_i(v_{i,k})\| \r)
     \sum_{j=1}^{m} \|y_k - w_{j,k}\|\cr &- 2 \alpha_{k+1}
     \sum_{i=1}^{m} \epsilon_{i,k+1}^T(v_{i,k} -z) + \alpha_{k+1}^2
     \sum_{i=1}^{m} \l\|\nabla f_i(v_{i,k}) + \epsilon_{i,k+1} \r\|^2.
   \end{align*}
   The result follows by using the subgradient norm boundedness, 
$\|\nabla f_i(v_{i,k})\|\le C_i$ for all $k$ and $i$.
\end{proof}

\section{Convergence in mean}
\label{sec:generror}
Here, we study the behavior of the iterates 
generated by the algorithm, under the
assumption that the errors have bounded norms in mean square.
In particular, we assume the following.
\begin{assumption}
  \label{ass:erindep1}
  The subgradient errors are uniformly bounded in mean square, i.e, there
  are scalars $\bar{\nu}_i$ such that
 \[ \EXP{\|\epsilon_{i,k+1}\|^2}
  \le \bar{\nu}_i^2\qquad\hbox{for all $i \in V$ and all $k$}.\]
\end{assumption}
Using this assumption, we provide a bound on the expected disagreement
$\EXP{\|w_{i,k}- y_{k}\|}$
for nondiminishing stepsize. We later use this bound to provide an
estimate for the algorithm's performance in mean. 
The bound is provided in the following theorem.

\begin{theorem}
  \label{thm:expconsensus}
  Let Assumptions~\ref{ass:convex}a, \ref{ass:con}, \ref{ass:weight}
  and \ref{ass:erindep1} hold.  Also, let the subgradients of each
  $f_i$ be uniformly bounded over $X$, i.e., for each $i\in V$ there
  is $C_i$ such that
  \[\|\nabla f_i(x)\|\le C_i\qquad\hbox{for all }x\in X.\] 
   If the stepsize $\{\alpha_k\}$ is such that
    $\lim_{k\to\infty} \alpha_{k}=\alpha$ for some $\alpha \geq 0,$ then 
    for all $j\in V$,
    \begin{align*}
      \limsup_{k\to\infty} \EXP{\|y_{k+1} - w_{j,k+1}\|} \leq
      \alpha\max_{i\in V}\{C_i+\bar{\nu}_i\}\,
      \l(2+\frac{m\theta\beta}{1-\beta}\r).
    \end{align*}
\end{theorem}
\begin{proof}
  The conditions of Lemma \ref{lemma:yandwj} are satisfied.  Taking
  the expectation in the relation of Lemma \ref{lemma:yandwj} 
  and using the inequality $\EXP{\| \epsilon_{i,k}\|}
  \leq \sqrt{\EXP{\| \epsilon_{i,k}\|^2}} = \bar{\nu}_i,$ we obtain
  for all $j\in V$ and all $k,$
  \begin{align}
    \EXP{\|y_{k+1} - w_{j,k+1}\|} \leq& m \theta \beta^{k+1}
    \max_{i\in V} \l\|w_{i,0}\r\| + m\theta\beta \max_{i\in V}\,
    \{C_i+\bar{\nu}_i\} \sum_{\ell=1}^{k} \beta^{k-\ell} \alpha_{\ell}
    \nn \\ &+2\alpha_{k+1}\max_{i\in V} \{C_i+\bar{\nu}_i\}. 
    %\leq& m \theta \beta^{k+1} \max_{i\in V} \l\|w_{i,0}\r\| +
    %m\theta\beta \max_{i\in V}\, \{C_i+\bar{\nu}_i\} \sum_{\ell=1}^{k}
    %\beta^{k-\ell} \l(\alpha_{\ell} - \alpha \r) \nn \\ & +
    %\theta\beta m \alpha \max_{i\in V}\, \{C_i+\bar{\nu}_i\}
    %\sum_{\ell=1}^{k} \beta^{k-\ell} + 2\alpha_{k+1}\max_{i\in V}
    %\{C_i+\bar{\nu}_i\}.
    \label{eqn:keyrel}
  \end{align}
  Since $\lim_{k\to\infty}\a_{k}=\alpha$, by
  Lemma \ref{lemma:betalpha}(a) we have
  $\lim_{k\to\infty}\sum_{\ell=1}^{k} 
  \beta^{k-\ell} \alpha_{\ell}=\frac{\alpha}{1-\beta}.$
  Using this relation and $\lim_{k\to\infty}\a_{k}=\alpha$, 
  we obtain the result by taking the limit superior in (\ref{eqn:keyrel}) 
  as $k\to\infty$.
\end{proof}

When the stepsize is diminishing (i.e., $\alpha = 0$),
the result of Theorem \ref{thm:expconsensus} implies that 
the expected disagreements 
$\EXP{\|y_{k+1} - w_{j,k+1}\|}$ converge to 0 
for all $j$. Thus, there is an asymptotic consensus in mean. 
We formally state this as a corollary.

\begin{corollary}
  \label{cor:expconsensus}
  Let the conditions of Theorem~\ref{thm:expconsensus} hold with
  $\alpha =0.$ Then $\lim_{k\to\infty} 
  \EXP{\|w_{j,k} - y_{k} \|}= 0$ for all $j \in V.$
\end{corollary}

We next obtain bounds on the performance of the algorithm. We
make the additional assumption that the set $X$ is bounded. Thus, the
subgradients of each $f_i$ are also bounded (see \cite{Bertsekas03},
Proposition~4.2.3). 

Note that, under Assumption~\ref{ass:erindep1},
by Jensen's inequality we have  
$\|\EXP{\epsilon_{i,k+1}}\|\le \bar{\nu}_i.$
Therefore, under Assumption~\ref{ass:erindep1},
\begin{equation}
\limsup_{k\to \infty}\|\EXP{\epsilon_{i,k+1}}\|\le \bar{\nu}_i
\qquad\hbox{for all }i\in V.
\label{eqn:boundmu}
\end{equation}
%$\bar{\mu}_{i,k} \leq \bar{\nu}_i.$
We have used this relation in our analysis of the agent disagreements 
in Theorem~\ref{thm:expconsensus}. 
%We will also use later on in
%the analysis of the consensus.  
%However, when studying the behavior of the algorithm we
%use $\limsup_{k\to \infty}\|\EXP{\epsilon_{i,k+1}}\|$ explicitly 
%to capture more refined impacts of the errors.  
Using this relation, we obtain special results for
the cases when the errors are zero mean or when their mean is
diminishing, i.e., the cases $\EXP{\epsilon_{i,k+1}}=0$ for all $i,k,$
or $\limsup_{k\to \infty}\|\EXP{\epsilon_{i,k+1}}\|=0$ for all $i$.

%%%%%%%%%%%%%%%%%%%%%%%%%%%%%%%%%%%%%%%%%%%%%%%%%%%%%%%%%%%%%%%%%%%%%%%%%%%%%%
\begin{theorem}
  \label{thm:constantstep}
  Let Assumptions~\ref{ass:convex}, \ref{ass:con}, \ref{ass:weight}
  and \ref{ass:erindep1} hold. Assume that the set $X$ is bounded.
  Let $\lim_{k\to\infty}\alpha_{k} = \alpha$ with $\alpha \geq 0.$ If
  $\alpha = 0,$ also assume that $\sum_{k} \alpha_k = \infty.$ Then,
  for all $j\in V$,
  \begin{eqnarray*}
    \liminf_{k \to \infty} \EXP{f(w_{j,k})} \leq f^* + \max_{x,y \in
      X}\l\|x - y\r\| \sum_{i=1}^{m} \bar{\mu}_i + m\alpha \l(\max_{i \in V}
      \{C_i + \bar{\nu}_i\} \r)^2 \l( \frac{9}{2} + \frac{2m \theta\beta} {1
      - \beta}\r),
  \end{eqnarray*}
  where $\bar\mu_i=\limsup_{k\to \infty}\|\EXP{\epsilon_{i,k+1}}\|$ and 
  $C_i$ is an upper-bound on the subgradient norms of $f_i$ over the
  set $X.$
\end{theorem}
\begin{proof}
  Under Assumption \ref{ass:erindep1}, the limit superiors 
  $\bar\mu_i=\limsup_{k\to \infty}\|\EXP{\epsilon_{i,k+1}}\|$ are finite 
  [cf.~Eq.~(\ref{eqn:boundmu})].
  Since the set $X$ is bounded the subgradients of $f_i$ over the set
  $X$ are also bounded for each $i\in V$; hence, the bounds $C_i,\
  i\in V$ on subgradient norms exist.  Thus, the conditions of
  Lemma~\ref{lemma:biralpha} are satisfied.  Further, by 
  Assumption~\ref{ass:convex}, the set $X$ is contained in the  
  interior of the domain of $f$, over which the function is continuous 
 (by convexity; see \cite{Rockafellar70}). 
  Thus, the set $X$ is compact and $f$ is continuous over $X$, implying that 
  the optimal set $X^*$ is nonempty. Let $x^* \in
  X^*,$ and let $y = x^*$ in Lemma~\ref{lemma:biralpha}. We have, for
  all $k$,
  \begin{eqnarray*}
    \sum_{i=1}^{m} \|v_{i,k+1} - x^*\|^2 &\leq& \sum_{i=1}^{m} \|
    v_{i,k} - x^*\|^2 -2 \alpha_{k+1}\, \l( f(y_{k}) - f^* \r) \cr
    &&+2\alpha_{k+1}\l(\max_{i\in V} C_i\r)\sum_{j=1}^m
    \|y_k-w_{j,k}\|\cr 
    &&-2\alpha_{k+1} \sum_{i=1}^{m}
    \epsilon_{i,k+1}^T(v_{i,k}-x^*) + \alpha_{k+1}^2 \sum_{i=1}^{m} \l(
    C_i + \|\epsilon_{i,k+1}\| \r)^2.
  \end{eqnarray*}
  Since $X$ is bounded, by using $\|v_{i,k}-x^*\|\le \max_{x,y\in
    X}\|x-y\|,$ taking the expectation and using the error bounds
    $\EXP{\|\epsilon_{i,k+1}\|^2}\le \bar \nu_i^2$ %and 
    %$\bar\mu_{i,k+1}=\l\|\EXP{\epsilon_{i,k+1}}\r\|,$
    we obtain
  \begin{eqnarray}
    \sum_{i=1}^{m} \EXP{\|v_{i,k+1} - x^*\|^2} &\leq& \sum_{i=1}^{m}
    \EXP{\| v_{i,k} - x^*\|^2} -2 \alpha_{k+1}\, \l( \EXP{f(y_{k})} -
    f^* \r) \cr 
    &&+2\alpha_{k+1}\l(\max_{i\in V} C_i\r)\sum_{j=1}^m
    \EXP{\|y_k-w_{j,k}\|}\cr 
    &&+2\alpha_{k+1} \,\max_{x,y\in X}\|x-y\|
    \sum_{i=1}^{m} \l\|\EXP{\epsilon_{i,k+1}}\r\|%\bar{\mu}_{i,k+1} 
    + \alpha_{k+1}^2 \sum_{i=1}^m \l(C_i + \bar{\nu}_i\r)^2. \qquad  
    \label{eqn:estim1} %\\
    %&\leq& \sum_{i=1}^{m}
    %  \EXP{\| v_{i,k} - x^*\|^2} \nn \\&&- 2\alpha_{k+1}\l( \l(\EXP{f(y_{k})}
    %  - f^*\r) -\l(\max_{i \in V} C_i \r) \sum_{j=1}^{m} \EXP{\|y_k
    %  -w_{j,k}\|}\r. \nn \\&&-\max_{x,y \in X}\|x -y\| \sum_{i=1}^{m}
    %  \bar{\mu}_{i,k+1} - \l. \frac{m\alpha_{k+1}}{2}\, \l(\max_{i\in V}\{
    %  C_i + \bar{\nu}_i\}\r)^2\r). \nn
   \end{eqnarray}
  By rearranging the terms and summing over $k=1,\ldots,K,$ for an 
  arbitrary $K$, we obtain
  \begin{align*}
      &2\sum_{k=1}^{K} \alpha_{k+1} \l( \l(\EXP{f(y_{k})} - f^*\r)
     -\l(\max_{i \in V} C_i \r) \sum_{j=1}^{m} \EXP{\|y_k
     -w_{j,k}\|}\r.\cr 
     & -\l. \max_{x,y \in X}\|x -y\|
     \sum_{i=1}^{m} \l\|\EXP{\epsilon_{i,k+1}}\r\|
     %\bar{\mu}_{i,k+1} 
     - \frac{m\alpha_{k+1}}{2}\,
     \l(\max_{i\in V}\{ C_i + \bar{\nu}_i\}\r)^2\r)\cr &\leq
     \sum_{i=1}^{m} \EXP{\| v_{i,1} - x^*\|^2} - \sum_{i=1}^{m}
     \EXP{\|v_{i,K+1} - x^*\|^2} \leq m \max_{x,y \in X} \|x-y\|^2.
    \end{align*}
    Note that when $\alpha_{k+1} \to \alpha$ and $\alpha > 0,$ we have
    $\sum_{k} \alpha_{k} = \infty.$ When $\alpha = 0,$ we have assumed
    that $\sum_{k} \alpha_{k} = \infty.$ 
    Therefore, by letting $K\to\infty$, we have
    \begin{align*}
      &\liminf_{k\to \infty} \l(\EXP{f(y_{k})} -\l(\max_{i \in V} C_i
      \r) \sum_{j=1}^{m} \EXP{\|y_k -w_{j,k}\|}\r.\cr 
      &\l.  -\max_{x,y
      \in X}\|x -y\| \sum_{i=1}^{m} \l\|\EXP{\epsilon_{i,k+1}}\r\|
       %\bar{\mu}_{i,k+1}
      -\frac{m\alpha_{k+1}}{2}\, \l(\max_{i\in V}\{ C_i +
      \bar{\nu}_i\}\r)^2\r) \leq f^*.
    \end{align*}
    Using $\limsup_{k\to \infty} \l\|\EXP{\epsilon_{i,k+1}}\r\|
    = \bar{\mu}_i$ [see Eq.~(\ref{eqn:boundmu})]
    and $\lim_{k\to \infty} \alpha_{k} = \alpha,$ we obtain
    \begin{align*}
      \liminf_{k\to \infty} \EXP{f(y_{k})} \le & f^* +
      \frac{m\alpha}{2}\, \l(\max_{i\in V}\{ C_i +
      \bar{\nu}_i\}\r)^2 +\l(\max_{i \in V} C_i \r) \sum_{j=1}^{m}
      \limsup_{k\to\infty} \EXP{\|y_k -w_{j,k}\|} \cr
      &+\max_{x,y \in X}\|x
      -y\| \sum_{i=1}^{m} \bar{\mu}_{i}.
    \end{align*}
     Next from the convexity inequality in (\ref{eqn:lip3}) and the
     boundedness of the subgradients it follows that for all $k$ and
     $j\in V$,
    \begin{align*}
      \EXP{f(w_{j,k}) - f(y_{k})} \leq \l(\sum_{i=1}^{m}C_i \r)
      \EXP{\|y_{k} - w_{j,k} \|},
    \end{align*}
    implying 
    \begin{align*}
      \liminf_{k\to \infty} \EXP{f(w_{j,k})} \le f^* & +
      \frac{m\alpha}{2}\, \l(\max_{i\in V}\{ C_i +
      \bar{\nu}_i\}\r)^2 +\l(\max_{i \in V} C_i \r) \sum_{j=1}^{m}
      \limsup_{k\to\infty} \EXP{\|y_k -w_{j,k}\|} \\ 
      &+
      \l(\sum_{i=1}^{m}C_i \r) \limsup_{k\to\infty} \EXP{\|y_{k} -
      w_{j,k} \|}  +\max_{x,y \in X}\|x -y\| \sum_{i=1}^{m}
      \bar{\mu}_{i}.
    \end{align*}
    By Theorem \ref{thm:expconsensus}, we have for all $j\in V$,
    \begin{align*}
      \limsup_{k\to\infty} \EXP{\|y_{k} - w_{j,k}\|} \leq
      \alpha\max_{i\in V}\{C_i+\bar{\nu}_i\}\,
      \l(2+\frac{m\theta\beta}{1-\beta}\r).
      %\label{eqn:limsupyk}
    \end{align*}
    By using the preceding relation, we see that
    \begin{align*}
      \liminf_{k\to \infty} \EXP{f(w_{j,k})} \leq& f^* +
      \frac{m\alpha}{2}\, \l(\max_{i\in V}\{ C_i +
      \bar{\nu}_i\}\r)^2 + \max_{x,y \in X}\|x -y\| \sum_{i=1}^{m}
      \bar{\mu}_{i} \\ &+m\alpha \l(\max_{i \in V} C_i \r) \max_{i\in
      V}\{C_i+\bar{\nu}_i\}\, \l(2+\frac{m\theta\beta}{1-\beta}\r)
      \\&+ \alpha \l(\sum_{i=1}^{m}C_i \r) \max_{j\in
      V}\{C_j+\bar{\nu}_j\}\, \l(2+\frac{m\theta\beta}{1-\beta}\r) \\
      &\hspace{- 0.15 in}\leq f^* + \max_{x,y \in X}\l\|x - y\r\|
      \sum_{i=1}^{m} \bar{\mu}_i + m\alpha 
      \l(\max_{i \in V} \{C_i +
      \bar{\nu}_i\}\r)^2 \l( \frac{9}{2} + \frac{2m \theta\beta} {1 -
      \beta}\r).
    \end{align*}
\end{proof}

The network topology influences
the error only through the term $ \frac{\theta\beta} {1 - \beta}$ and
can hence be used as a figure of merit for comparing different 
topologies. For a network that is strongly connected at every time,
[i.e., $Q=1$ in
Assumption~\ref{ass:con}] and when $\eta$ in Assumption
\ref{ass:weight} does not depend on the number $m$ of agents, the term
$\frac{\theta\beta} {1 - \beta}$ is of the order $m^2$ and the error
bound scales as $m^4.$

We next show that stronger bounds can be obtained for a specific
weighted time averages of the iterates $w_{i,k}$. In particular, we
investigate the limiting behavior of $\{f(z_{i,t})\},$ where ${z_{i,t}
= }\frac{\sum_{k=1}^{t} \alpha_{k+1} w_{i,k}}{\sum_{k=1}^{t}
\alpha_{k+1}}.$ Note that agent $i$ can locally and recursively
evaluate $z_{i,t+1}$ from $z_{i,t}$ and $w_{i,t+1}.$

\begin{theorem}
  \label{thm:averi}
Consider the weighted time averages ${z_{j,t}
= }\frac{\sum_{k=1}^{t} \alpha_{k+1} w_{j,k}}{\sum_{k=1}^{t}
\alpha_{k+1}}$ for $j\in V$ and $t\ge1.$ %of the iterates $w_{j,k}$.  
Let the conditions of Theorem~\ref{thm:constantstep} hold.
Then, we have for all $j\in V,$
  \[
  \limsup_{t\to\infty} \EXP{f\l(z_{j,t}\r)} \leq f^* + \max_{x,y \in
  X} \|x-y\| \sum_{i=1}^{m} \bar{\mu}_{i} + m \alpha \l(\max_{i\in V}
  \{C_i + \bar{\nu}_i\}\r)^2 \left(\frac{9}{2}
  +\frac{2m\theta\beta}{1-\beta}\r).
  \]
%  with $\bar\mu_i= \limsup_{k\to\infty}\|\EXP{\epsilon_{i,k+1}}\|$. 
\end{theorem}
%\selfnote{Changed the proof. Also added a Lemma 3.1c}

\begin{proof}
The relation in (\ref{eqn:estim1}) of Theorem~\ref{thm:constantstep}
is valid, and we have for any $x^*\in X^*,$
\begin{eqnarray*}
  \sum_{i=1}^{m} \EXP{\|v_{i,k+1} - x^*\|^2} &\leq& \sum_{i=1}^{m}
  \EXP{\| v_{i,k} - x^*\|^2} -2 \alpha_{k+1}\, \l( \EXP{f(y_{k})} -
  f^* \r) \\ &&+2\alpha_{k+1}\l(\max_{i\in V} C_i\r)\sum_{\ell=1}^m
  \EXP{\|y_k-w_{\ell,k}\|} \\ &&+2\alpha_{k+1}\max_{x,y \in X} \|x-y\|
  \sum_{i=1}^{m} \|\EXP{\epsilon_{i,k+1}}\|
  + \alpha_{k+1}^2 \sum_{i=1}^{m}(C_i + \bar{\nu}_i)^2.
\end{eqnarray*}
From the subgradient boundedness and the
subgradient inequality in (\ref{eqn:lip3}) we have for any $j$,
\[
\EXP{f(y_{k})} - \EXP{f(w_{j,k})} \geq -\l(\sum_{i=1}^{m} C_i \r)
\EXP{\|y_{k}-w_{j,k}\|} \geq - m \l( \max_{i\in V} C_i \r)
\EXP{\|y_{k}-w_{j,k}\|}.
\]
Therefore, we obtain 
\begin{eqnarray*}
  \sum_{i=1}^{m} \EXP{\|v_{i,k+1} - x^*\|^2} &\leq& \sum_{i=1}^{m}
  \EXP{\| v_{i,k} - x^*\|^2} -2 \alpha_{k+1}\, \l( \EXP{f(w_{j,k})} -
  f^* \r) \\ &&+2\alpha_{k+1}\l(\max_{i\in V}
  C_i\r)\l(m\EXP{\|y_k-w_{j,k}\|} + \sum_{i=1}^m
  \EXP{\|y_k-w_{i,k}\|}\r) \\ &&+2\alpha_{k+1}\max_{x,y \in X} \|x-y\|
  \sum_{i=1}^{m} \|\EXP{\epsilon_{i,k+1}}\|
  + \alpha_{k+1}^2 \sum_{i=1}^{m}(C_i + \bar{\nu}_i)^2.
\end{eqnarray*}
By re-arranging these terms, summing over $k=1,\ldots,t$ and dividing
with $2\sum_{k=1}^{t} \alpha_{k+1}$, we further obtain
{\small
  \begin{eqnarray*}
    \sum_{k=1}^t \frac{\alpha_{k+1} \EXP{f(w_{j,k})}}{\sum_{k=1}^{t}
      \alpha_{k+1}}\, &\leq& f^* + \frac{1}{2\sum_{k=1}^{t}
      \alpha_{k+1}}\sum_{i=1}^{m} \EXP{\| v_{i,1} - x^*\|^2} \\
      &&+\sum_{k=1}^{t} \frac{\alpha_{k+1}\l(\max_{i\in V}
      C_i\r)\l(m\EXP{\|y_k-w_{j,k}\|} + \sum_{i=1}^m
      \EXP{\|y_k-w_{i,k}\|}\r)}{\sum_{k=1}^{t} \alpha_{k+1}} \\ &+&
      \max_{x,y \in X} \|x-y\| 
      \sum_{i=1}^{m} \frac{\sum_{k=1}^{t}
      \alpha_{k+1} \|\EXP{\epsilon_{i,k+1}}\|}{\sum_{k=1}^{t} \alpha_{k+1}} +
      \frac{\sum_{k=1}^{t} \alpha_{k+1}^2}{2\sum_{k=1}^{t}
      \alpha_{k+1}} \sum_{i=1}^{m} (C_i + \bar{\nu}_i)^2.
\end{eqnarray*}}
Next by the convexity of $f$ note that
\[
f(z_{j,t}) = f\l(\sum_{k=1}^t \frac{\alpha_{k+1}
  w_{j,k}}{\sum_{k=1}^{t} \alpha_{k+1}}\r) \leq \sum_{k=1}^t
  \frac{\alpha_{k+1} f(w_{j,k})}{\sum_{k=1}^{t} \alpha_{k+1}}.
\]
From the preceding two relations we obtain
{\small
  \begin{eqnarray}
  \EXP{f(z_{j,t})} &\leq& f^* + \frac{1}{2\sum_{k=1}^{t}
   \alpha_{k+1}}\sum_{i=1}^{m} \EXP{\| v_{i,1} - x^*\|^2} \nn \\
   &&+\sum_{k=1}^{t} \frac{\alpha_{k+1}\l(\max_{i\in V}
   C_i\r)\l(m\EXP{\|y_k-w_{j,k}\|} + \sum_{i=1}^m
   \EXP{\|y_k-w_{i,k}\|}\r)}{\sum_{k=1}^{t} \alpha_{k+1}}\nn \\ &&+
   \max_{x,y \in X} \|x-y\| \sum_{i=1}^{m} \frac{\sum_{k=1}^{t}
   \alpha_{k+1} \|\EXP{\epsilon_{i,k+1}}\|}{\sum_{k=1}^{t} \alpha_{k+1}}  +
   \frac{\sum_{k=1}^{t} \alpha_{k+1}^2}{2\sum_{k=1}^{t} \alpha_{k+1}}
   \sum_{i=1}^{m} (C_i + \bar{\nu}_i)^2.
   \qquad\quad\label{eqn:lim0}
\end{eqnarray}}
\hskip -0.4pc First note that in the limit as $t \to \infty,$ 
the second term in (\ref{eqn:lim0}) converges to $0$ 
since $\sum_{k=1}^{t}\alpha_{k+1} = \infty.$ 
By using the results of Lemma~\ref{lemma:betalpha}c, 
for the remaining terms, we obtain
\begin{eqnarray*}
  \limsup_{t \to \infty} \EXP{f(z_{j,t})} &\leq& f^* 
   +\l(\max_{i\in V} C_i\r) 
   \limsup_{k\to\infty}\l(m\EXP{\|y_k-w_{j,k}\|} + \sum_{i=1}^m
   \EXP{\|y_k-w_{i,k}\|}\r)
   \nn\\ &&+\max_{x,y \in X} \|x-y\|\sum_{i=1}^{m}
   \limsup_{k\to\infty}\|\EXP{\epsilon_{i,k+1}}\|
   + \frac{\alpha}{2}
   \sum_{i=1}^{m} (C_i + \bar{\nu}_i)^2.
\end{eqnarray*}
By Theorem \ref{thm:expconsensus}, we have for all $j\in V$,
    \begin{align*}
      \limsup_{k\to\infty} \EXP{\|y_{k} - w_{j,k}\|} \leq
      \alpha\max_{i\in V}\{C_i+\bar{\nu}_i\}\,
      \l(2+\frac{m\theta\beta}{1-\beta}\r),
      %\label{eqn:limsupyk}
    \end{align*} 
which when substituted in the preceding relation, yields
\begin{eqnarray*}
  \limsup_{t \to \infty} \EXP{f(z_{j,t})} &\leq& f^* +
   2m\alpha\l(\max_{i\in V} C_i\r) \max_{i\in V} \{C_i+\bar{\nu}_i\}
   \l(2+\frac{m\theta\beta}{1 - \beta}\r) \nn\\ 
   &&+\max_{x,y \in X} \|x-y\| \sum_{i=1}^{m}
   \limsup_{k\to\infty}\|\EXP{\epsilon_{i,k+1}}\|
   + \frac{\alpha}{2}\sum_{i=1}^{m} (C_i + \bar{\nu}_i)^2 \\ 
   &\leq& f^* + \max_{x,y \in X} \|x-y\| \sum_{i=1}^{m}
   \limsup_{k\to\infty}\|\EXP{\epsilon_{i,k+1}}\|\\
   &&+ m \alpha \l(\max_{i\in V}
   \{C_i + \bar{\nu}_i\}\r)^2 \left(\frac{9}{2}
   +\frac{2m\theta\beta}{1-\beta}\r).
\end{eqnarray*}
\end{proof}

The error bounds in Theorems~\ref{thm:constantstep} and
\ref{thm:averi} have the same form, but they apply 
to different sequences of function evaluations. Furthermore, in  
Theorem~\ref{thm:averi}, the bound is for \emph{all}
subsequences of $\EXP{f(z_{i,k})}$ for each agent $i$. In
contrast, in Theorem~\ref{thm:constantstep}, the bound
is only for \emph{a} subsequence of $\EXP{f(z_{i,k})}$ for each agent~$i$. 
Theorem~\ref{thm:averi} demonstrates that, due to the convexity
of the objective function $f$, there is an advantage 
when agents are using the running averages of their iterates.

When the error\footnote{
When the moments $\|\EXP{\epsilon_{i,k+1}}\|$ are zero, 
it can be seen that the results of Theorems~\ref{thm:constantstep} 
and \ref{thm:averi}
hold when the boundedness of $X$
is replaced by the weaker assumption that the subgradients of each
$f_i$ are bounded over $X.$}
moments $\|\EXP{\epsilon_{i,k+1}}\|$ converge 
to zero as $k\to\infty$,
and the stepsize converges to zero [$\a=0$],
Theorems~\ref{thm:constantstep} and \ref{thm:averi} yield respectively 
\[\liminf_{k\to\infty}\EXP{f(w_{j,k})} = f^*\qquad\hbox{and}\qquad
\lim_{k\to\infty}\EXP{f(z_{j,k})} = f^*.\] 

When a constant stepsize $\alpha$ is used, the vector 
$z_{j,t}$ is simply the
running average of all the iterates of agent $j$ until time $t,$ i.e., 
$z_{j,t}= \frac{1}{t} \sum_{k=1}^{t} w_{j,k}.$ For this case, with zero mean
errors, the relation in (\ref{eqn:lim0}) reduces to
\begin{eqnarray}
  \EXP{f(z_{j,t})} &\leq& f^* + \frac{1}{2t\alpha} \sum_{i=1}^{m}
   \EXP{\| v_{i,1} - x^*\|^2}\nn \\ &&+ \l(\max_{i\in V} C_i\r)
   \frac{1}{t} \sum_{k=1}^{t} \l(m\EXP{\|y_k-w_{j,k}\|} + \sum_{i=1}^m
   \EXP{\|y_k-w_{i,k}\|}\r)\nn \\ && + \frac{\alpha}{2} 
   \sum_{i=1}^{m} (C_i +
   \bar{\nu}_i)^2. \label{eqn:cs}
\end{eqnarray}
This can be used to derive an estimate per iteration, as 
seen in the following.

\begin{corollary}\label{cor:boundstep}
Under the conditions of Theorem~\ref{thm:constantstep} with
$\|\EXP{\epsilon_{i,k+1}}\|=0$ and $\a_k=0$ for all $i$ and $k$,
for the average sequences $\{z_{j,k}\}$ we have for all $t$ and $j$,
\begin{eqnarray*}
  \EXP{f(z_{j,t})} &\leq& f^* + \frac{1}{2t\alpha} \sum_{i=1}^{m}
   \EXP{\| v_{i,1} - x^*\|^2} 
   + \frac{2m^2\theta \beta^2}{t(1 -
   \beta)}\l(\max_{i \in V} C_i \r)\l(\max_{i \in V} \|w_{i,0}\|\r) \nn \\
   &&+ m\alpha \l(\max_{i\in V}\{C_i + \bar{\nu}_i\}\r)^2
   \l(\frac{9}{2} + \frac{ 2m\theta
   \beta}{1 - \beta} \r).
\end{eqnarray*}
\end{corollary}

\begin{proof}
Taking the expectation in the relation of Lemma~\ref{lemma:yandwj}, we obtain
\begin{align*}
    \EXP{\|y_{k+1} - w_{j,k+1}\|} \leq& m \theta \beta^{k+1}
    \max_{i\in V} \l\|w_{i,0}\r\| + m \alpha \theta\beta 
\l(\max_{i\in
    V}\, \{C_i+\bar{\nu}_i\}\r) \sum_{\ell=1}^{k} \beta^{k-\ell} \nn \\
    &+ 2\alpha \max_{i\in V} \{C_i+\bar{\nu}_i\} \\
    \leq& m \theta
    \beta^{k+1} \max_{i\in V} \l\|w_{i,0}\r\| 
    + \a \l(\max_{i\in V} \{C_i+\bar{\nu}_i\}\r)\left(2+\frac{m
    \theta\beta}{1 - \beta}\right).
\end{align*}
Combining the preceding relation 
with the inequality in (\ref{eqn:cs}), 
and using $\sum_{k=1}^t\beta^{k+1}\le \frac{\beta^2}{1-\beta}$,
we obtain 
\begin{eqnarray*}
  \EXP{f(z_{j,t})} &\leq& 
    f^* + \frac{1}{2t\alpha} \sum_{i=1}^{m}
   \EXP{\| v_{i,1} - x^*\|^2} + \frac{2m^2\theta \beta^2}{t(1 -
   \beta)}\l(\max_{i \in V} C_i \r)\l(\max_{i \in V} \|w_{i,0}\|\r) \nn \\
   &&+ m\alpha \l(\max_{i\in V}\{C_i + \bar{\nu}_i\}\r)^2
   \l(\frac{9}{2} + \frac{ 2m\theta
   \beta}{1 - \beta} \r).  %\label{eqn:fs}
\end{eqnarray*}
\end{proof}

The preceding equation provides a bound on the
algorithm's performance at each iteration. The bound
can be used in obtaining stopping rules for the algorithm. For
example, consider the error free case ($\bar{\nu}_i = 0$) and suppose
that the goal is to determine the number of iterations required 
for agents  to find a point in the $\epsilon$-optimal set,
i.e., in the set $X_{\epsilon} = \{x \in X: f(x) \le f^* + \epsilon\}.$
Minimizing the bound in Corollary~\ref{cor:boundstep} 
over different stepsize values
$\alpha$, we can show that $\epsilon$-optimality can be achieved in
$N_\epsilon=\l\lceil \frac{1}{\psi^2_{\epsilon}} \r\rceil$ iterations with a
stepsize $\alpha_{\epsilon} =
\frac{\sqrt{A}\psi_{\epsilon}}{\sqrt{C}}$, where $\psi_{\epsilon}$ is
the positive root of the quadratic equation
\[
Bx^2 + 2\sqrt{AC}x - \epsilon = 0,
\]
and $A, B$ and $C$ are
\begin{align*}
  &A = \frac{1}{2} \sum_{i=1}^{m}\| v_{i,1} - x^*\|^2,
  \qquad 
  B =
  \frac{2m^2\theta \beta^2}{1 - \beta}\l(\max_{i \in V} C_i
  \r)\l(\max_{i \in V} \|w_{i,0}\|\r),\\ 
  &C = m \l(\max_{i\in V}\{C_i +
  \bar{\nu}_i\}\r)^2\l(\frac{9}{2} + \frac{ 2m\theta \beta}{1 - \beta} \r).
\end{align*}
Since $\psi_{\epsilon}$ scales as $\sqrt{\epsilon},$ we can conclude
that $N_{\epsilon}$ scales as $\frac{1}{\epsilon^2}.$ Equivalently, we
can say that the level $\epsilon$ 
of sub-optimality diminishes inversely with
the square root of the number of iterations.

\section{Almost sure and mean square convergence}
\label{sec:dimerror} 
In this section, we impose some
additional assumptions on the subgradient errors to  
obtain almost sure consensus among the agents and 
almost sure convergence of the iterates to an optimal solution 
of (\ref{eqn:problem}).
Towards this, define $F_{k}$ to be the $\sigma$-algebra
$\sigma\l(\epsilon_{i,\ell};\ i\in V, 0 \leq \ell\leq k \r)$ generated
by the errors in the agent system up to time $k$. In other words,
$F_{k}$ captures the history of the errors until the end of
time $k$. We use the following assumption on the subgradient 
errors $\epsilon_{i,k}$.

\begin{assumption}
 \label{ass:erindep}
  %For all $i\in V,$ we have
 %$\sum_{k=0}^{\infty} \|\EXP{\epsilon_{i,k+1} \mid
    %F_{k}}\|^2 < \infty$ with probability 1. 
  There are
    scalars $\nu_i$ such that $ \EXP{\|\epsilon_{i,k+1}\|^2 \mid
    F_{k}} \le \nu_i^2$ for all $k$ with probability 1.
\end{assumption}

%Imposing assumptions on the conditional
%moments, rather than the absolute moments, is quite standard in
%stochastic optimization literature \cite{Bertsekas00}.  
Note that
Assumption~\ref{ass:erindep} is stronger than Assumption~\ref{ass:erindep1}. 
Furthermore, when the errors
are independent across iterations and across agents,
Assumption~\ref{ass:erindep} reduces to Assumption~\ref{ass:erindep1}. 

We start by analyzing the agents' disagreements measured in terms of distances
$\|y_k-w_{j,k}\|$. We have the following result.

\begin{theorem}
  \label{thm:asconsensus}
  Let Assumptions~\ref{ass:convex}a, \ref{ass:con}, 
  \ref{ass:weight} and~\ref{ass:erindep} hold.  
  Suppose that the subgradients of each $f_i$
  are uniformly bounded over $X$, i.e., for each $i\in V$ there is
  $C_i$ such that
  \[\|\nabla f_i(x)\|\le C_i\qquad\hbox{for all }x\in X.\] 
  %Also,
  %assume that, for each $i\in V$, there is a scalar $\nu_i$ such that 
  %\[
  %\EXP{\|\epsilon_{i,k+1}\|^2 \mid F_{k}} \le \nu^2_i \qquad\hbox{for
  %all $k$}.
  %\]
  If $\sum_{k=0}^{\infty} \alpha_{k+1}^2 <\infty,$
  then with probability 1, 
  $$ \sum_{k=1}^{\infty} \alpha_{k+2}\|y_{k+1} - w_{j,k+1}\| < \infty 
  \qquad\hbox{for all $j\in V$}.$$ 
  Furthermore, for all $j\in V$, 
  we have $ \lim_{k\to\infty} \|y_{k+1} -w_{j,k+1}\|=0$ 
  with probability 1 and in mean square.
\end{theorem}

\begin{proof} 
  By Lemma \ref{lemma:yandwj} and the subgradient boundedness, we have
  for all $j\in V$,
  \begin{align*}
    \|y_{k+1} - w_{j,k+1}\| \leq& m \theta \beta^{k+1} \max_{i\in V}
    \l\|w_{i,0}\r\| + \theta \sum_{\ell=1}^{k} \beta^{k+1-\ell}
    \sum_{i=1}^{m}\alpha_{\ell} \l(C_i + \|\epsilon_{i,\ell}\|\r)
    \cr &+ \frac{1}{m} \sum_{i=1}^{m} \alpha_{k+1}
    \l(C_i+\|\epsilon_{i,k+1}\| \r)+ \alpha_{k+1} \l(C_j
    +\|\epsilon_{j,k+1}\| \r).
  \end{align*}
  Using the inequalities $$\alpha_{k+2} \alpha_{\ell}
  \l(C_i+\|\epsilon_{i,\ell}\| \r) \leq \frac{1}{2}\l(\alpha^2_{k+2}
  +\alpha_{\ell}^2 \l(C_i+\|\epsilon_{i,\ell}\|\r)^2\r)$$ and
  $\l(C_i+\|\epsilon_{i,\ell}\|\r)^2\le 2C_i^2
  +2\|\epsilon_{i,\ell}\|^2,$ we obtain
  \begin{align*}
    \alpha_{k+2} \|y_{k+1} - w_{j,k+1}\| \leq& \alpha_{k+2}\, m \theta
    \beta^{k + 1} \max_{i\in V}\|w_{i,0}\| \cr &+
    \theta\sum_{\ell=1}^{k} \beta^{k+1 - \ell} \sum_{i=1}^{m} \l(
    \frac{1}{2} \alpha^2_{k+2} +\alpha_{\ell}^2
    \l(C_i^2+\|\epsilon_{i,\ell}\|^2\r) \r)\cr & + \frac{1}{m}
    \sum_{i=1}^{m} \l(\frac{1}{2} \alpha^2_{k+2} +\alpha_{k+1}^2
    \l(C_i^2+\|\epsilon_{i,k+1}\|^2\r) \r)\cr &+
    \frac{1}{2}\alpha^2_{k+2} + \alpha_{k+1}^2
    \l(C_j^2+\|\epsilon_{j,k+1}\|^2\r).
  \end{align*} 
  By using the inequalities $\sum_{\ell=1}^{k} \beta^{k+1 - \ell}\le
  \frac{\beta}{1-\beta}$ for all $k\ge1$ \ and \ $\frac{1}{2m}
  +\frac{1}{2}\le 1$, and by grouping the terms accordingly, from the
  preceding relation we have
  \begin{align*}
    \alpha_{k+2} \|y_{k+1} - w_{j,k+1}\| \leq& \alpha_{k+2} \, m
    \theta \beta^{k + 1} \max_{i \in V}\|w_{i,0}\| + \l( 1 + \frac{m
    \theta\beta }{2(1 -\beta)} \r) \alpha^2_{k+2}\cr & +\theta
    \sum_{\ell = 1}^{k}\alpha_{\ell}^2 \beta^{k+1 - \ell}
    \sum_{i=1}^{m} \l(C^2_i +\|\epsilon_{i,\ell} \|^2\r) \cr &
    +\frac{1}{m}\alpha_{k+1}^2 \sum_{i=1}^{m} \l(C_i^2+\|
    \epsilon_{i,k+1}\|^2 \r)+
    \alpha_{k+1}^2\l(C_j^2+\|\epsilon_{j,k+1} \|^2\r).
  \end{align*}
  Taking the conditional expectation and using
  $\EXP{\|\epsilon_{i,\ell}\|^2\mid F_{\ell-1}}\le \nu_i^2,$ and then taking
  the expectation again, we obtain
  \begin{align*}
    \EXP{\alpha_{k+2} \|y_{k+1} - w_{j,k+1}\|} \leq& \alpha_{k+2} \, m
    \theta \beta^{k + 1} \max_{i \in V}\|w_{i,0}\| + \l( 1 + \frac{m
    \theta\beta }{2(1 -\beta)} \r) \alpha^2_{k+2}\cr & +\theta \l(
    \sum_{i=1}^{m} \l(C^2_i +\nu_i^2\r) \r) \sum_{\ell = 1}^{k}
    \alpha_{\ell}^2 \,\beta^{k+1 - \ell}\cr & +
    \frac{1}{m}\alpha_{k+1}^2 \sum_{i=1}^{m} \l(C_i^2+\nu_i^2 \r)+
    \alpha_{k+1}^2\l(C_j^2+\nu_j^2\r).
    %\label{eqn:expalpu}
  \end{align*} 
  Since $\sum_{k} \alpha^2_{k} < \infty$ (and hence $\{\alpha_k\}$ bounded), 
  the first two terms and the last
  two terms are summable.  Furthermore, in view of Lemma
  \ref{lemma:betalpha} [part (b)], we have
  \[\sum_{k=1}^\infty
  \sum_{\ell = 1}^{k}\beta^{k+1 - \ell} \alpha_{\ell}^2<\infty.
  \]
  Thus, the third term is also summable. Hence $\sum_{k=1}^{\infty}
  \EXP{\alpha_{k+2} \|y_{k+1} - w_{j,k+1}\|} < \infty.$ 
  From the monotone
  convergence theorem \cite{Billingsley79}, it follows that
  \[\EXP{\sum_{k=1}^{\infty} \alpha_{k+2} \|y_{k+1} - w_{j,k+1}\|}=
\sum_{k=1}^{\infty} \EXP{ \alpha_{k+2} \|y_{k+1} - w_{j,k+1}\|},\] 
 and it is hence finite for all $j$. 
 If the expected value of a random variable
is finite, then the variable has to be finite with probability 1; thus,
with probability 1,
  \begin{equation}
    \sum_{k=1}^{\infty} \alpha_{k+2} \|y_{k+1} - w_{j,k+1}\| < \infty
\qquad\hbox{for all }j\in V.
    \label{eqn:onec}
  \end{equation}

  We now show that $\lim_{k\to\infty} \|y_{k} - w_{j,k}\|=0$ with
  probability 1 for all $j\in V.$ Note that the conditions of
  Theorem~\ref{thm:expconsensus} are satisfied with $\bar{\nu}_i =
  \nu_i$ and $\alpha = 0.$ Therefore, $\|y_{k} - w_{j,k}\|$ converges
  to $0$ in the mean and from (Fatou's) Lemma~\ref{lemma:Fat} it follows that
  \[
  0 \leq \EXP{\liminf_{k \to \infty} \l\| y_{k}-w_{j,k} \r\|
  } \leq \liminf_{k \to \infty} \EXP{\l\| y_{k} - w_{j,k}\r\|
  } = 0,
  \]
  and hence $ \EXP{\liminf_{k \to \infty}\l\| y_{k} - w_{j,k}\r\|}=0.$ 
  Therefore, with probability 1,
  \begin{align}
    \liminf_{k \to \infty} \l\| y_{k}-w_{j,k}\r\| =0. \label{eqn:dkey}
  \end{align}
  
  To complete the proof, in view of (\ref{eqn:dkey}) 
  it suffices to show that $\l\|y_{k}-w_{j,k} \r\|$
  converges with probability~1. To show this, we define
  \[
  r_{i,k+1} = \sum_{j=1}^{m} a_{i,j}(k+1) w_{j}(k) - \alpha_{k+1}
    \left(\nabla f_i\l(v_{i,k}\r) + \epsilon_{i,k+1} \right),
  \]
  and note that $P_{X}[r_{i,k+1}] = w_{i,k+1}$ 
  [see (\ref{eqn:sagrad}) and (\ref{eqn:averi})]. Since
  $y_k=\frac{1}{m}\sum_{i=1}^m w_{i,k}$ and the set $X$ is convex, it
  follows that $y_k\in X$ for all $k$. Therefore, by the non-expansive
  property of the Euclidean projection in~(\ref{eqn:proj1}), we have
  $\|w_{i,k+1} - y_{k}\|^2 \leq \| r_{i,k+1} - y_{k} \|^2$ for all
  $i\in V$ and all $k$. Summing these relations over all $i,$ we
  obtain
  \begin{align*}
    \sum_{i=1}^{m}\|w_{i,k+1} - y_{k}\|^2
    \leq \sum_{i=1}^{m} \| r_{i,k+1} - y_{k} \|^2
    \qquad\hbox{for all $k$.} 
  \end{align*}
  From $y_{k+1}=\frac{1}{m}\sum_{i=1}^m w_{i,k+1}$ and the fact that
  the average of vectors minimizes the sum of distances between each
  vector and arbitrary vector in $\Re^n$ [cf.~Eq~(\ref{eqn:ineq3})],
  we further obtain
  \[
  \sum_{i=1}^{m}\|w_{i,k+1} - y_{k+1}\|^2
  \le \sum_{i=1}^{m} \|w_{i,k+1} - y_{k}\|^2.
  \]
  Therefore, for all $k$,
  \begin{align}
    \sum_{i=1}^{m}\|w_{i,k+1} - y_{k+1} \|^2\leq 
    \sum_{i=1}^{m} \|r_{i,k+1} - y_{k} \|^2. 
  \label{eqn:dad2}
  \end{align}
  
  We next relate $\sum_{i=1}^{m} \|r_{i,k+1} - y_{k} \|^2$ to
  $\sum_{i=1}^{m} \|w_{i,k} - y_{k} \|^2.$ From the definition of
  $r_{i,k+1}$ and the equality $\sum_{j=1}^m a_{i,j}(k+1) = 1$
  [cf.~Assumption~\ref{ass:weight}b], we have
  \begin{align*}
    r_{i,k+1} - y_{k} =& \sum_{j=1}^m a_{i,j}(k+1) \l( w_{j,k} -
    y_{k}\r) - \alpha_{k+1} \left(\nabla f_{i} (v_{i,k}) +
    \epsilon_{i,k+1}\right)
 \end{align*}
  By Assumption~\ref{ass:weight}a and \ref{ass:weight}b, we have that
  the weights $a_{i,j}(k+1), j\in V$ yield a convex combination. Thus, 
  by the convexity of the norm [(\ref{eqn:ineq5}) and (\ref{eqn:ineq6})] and 
  by the subgradient boundedness,
  we have
  \begin{align*}
    \|r_{i,k+1} - y_{k} \|^2\leq& \sum_{j=1}^m a_{i,j}(k+1) \left\|
    w_{j,k} - y_{k}\right\|^2 + \alpha^2_{k+1} \left\|\nabla f_{i}
    (v_{i,k}) + \epsilon_{i,k+1}\right\|^2 \cr &+ 2 \alpha_{k+1}
    \left\|\nabla f_{i} (v_{i,k}) + \epsilon_{i,k+1}\right\|
    \sum_{j=1}^m a_{i,j}(k+1) \left\|w_{j,k} - y_{k}\right\|\cr &\le
    \sum_{j=1}^m a_{i,j}(k+1) \left\| w_{j,k} - y_{k}\right\|^2 +
    2\alpha^2_{k+1} \l(C_i^2+ \|\epsilon_{i,k+1}\|^2 \r)\cr &+ 2
    \alpha_{k+1}\l(C_i+\|\epsilon_{i,k+1}\|\r) \sum_{j=1}^m
    a_{i,j}(k+1) \left\|w_{j,k} - y_{k}\right\|.
  \end{align*}
  Summing over all $i$ and using $\sum_{i=1}^{m} a_{i,j}(k+1) = 1$
  [cf.~Assumption~\ref{ass:weight}d], we obtain
  \begin{align*}
    \sum_{i=1}^{m} \|r_{i,k+1} - y_{k} \|^2 \leq& \sum_{j=1}^m \left\|
    w_{j,k} - y_{k}\right\|^2 + 2\alpha^2_{k+1} \sum_{i=1}^{m}
    \l(C_i^2+\|\epsilon_{i,k+1}\|^2\r) \cr &+ 2 \alpha_{k+1}
    \sum_{i=1}^{m} \l(C_i+\|\epsilon_{i,k+1}\|\r) \sum_{j=1}^m
    a_{i,j}(k+1) \left\| w_{j,k} - y_{k} \right\|.
  \end{align*}  
  Using this in (\ref{eqn:dad2}) and taking the conditional expectation, 
  we see that for all $k$, we have with probability 1,  
  \begin{align}
    \sum_{i=1}^{m} \EXP{\|w_{i,k+1} - y_{k+1} \|^2\mid F_k} \leq&
    \sum_{i=1}^m \left\| w_{i,k} - y_{k}\right\|^2 + 2\alpha^2_{k+1}
    \sum_{i=1}^{m} \l(C_i^2+\nu_i^2\r) \cr &+ 2 \alpha_{k+1}
    \sum_{i=1}^{m} \l(C_i+\nu_i\r) \sum_{j=1}^m \left\| w_{j,k} -
    y_{k} \right\|,
    \label{eqn:dad5}
  \end{align}
  where we use $a_{i,j}(k+1) \leq 1$ for all $i,j$ and $k$, 
and the relations  $\EXP{\|\epsilon_{i,k+1}\|^2\mid F_k}\le\nu_i^2,$ 
  $\EXP{\|\epsilon_{i,k+1}\|\mid F_k}\le\nu_i$ 
holding with probability 1. 

  We now apply Theorem~\ref{lemma:asm} to the relation in
    (\ref{eqn:dad5}).  To verify that the conditions of
    Theorem~\ref{lemma:asm} are satisfied, note that the stepsize
    satisfies $\sum_{k=1}^{\infty} \alpha^2_{k+1}<\infty$ for all
    $i\in V.$ We also have $\sum_{k=1}^{\infty} \alpha_{k+1} \left\|
    w_{j,k} - y_{k} \right\| < \infty$ with probability 1
    [cf.~(\ref{eqn:onec})].  Therefore, the relation
    in~(\ref{eqn:dad5}) satisfies the conditions of
    Theorem~\ref{lemma:asm} with $\zeta_k = D_k = 0,$ thus implying
    that $\left\| w_{j,k} - y_{k} \right\|$ converges with probability
    1 for every $j\in V.$
\end{proof}

Let us compare Theorem~\ref{thm:asconsensus} and
Corollary~\ref{cor:expconsensus}. Corollary~\ref{cor:expconsensus}
provided sufficient conditions for the different agents to have
consensus in the mean. Theorem~\ref{thm:asconsensus} strengthens this
to consensus with probability 1 and in mean square sense, for a
smaller class of stepsize sequences under a stricter assumption. 

We next show that the consensus vector is actually in the optimal set,
provided that the optimal set is nonempty
and the 
conditional expectations $\|\EXP{\epsilon_{i,k+1} \mid F_{k}}\|$
are diminishing.

\begin{theorem}
  \label{thm:dss}
  Let Assumptions~\ref{ass:convex}, \ref{ass:con}, \ref{ass:weight}
  and \ref{ass:erindep} hold. Suppose that the subgradients of each $f_i$
  are uniformly bounded over $X$, i.e., for each $i\in V$ there is
  $C_i$ such that
  \[\|\nabla f_i(x)\|\le C_i\qquad\hbox{for all }x\in X.\]
  Also, assume that $\sum_{k=0}^\infty\|\EXP{\epsilon_{i,k+1} \mid F_{k}}\|^2<\infty$
  for all $i\in V.$
  Further, let the stepsize sequence $\{\alpha_k\}$  be 
  such that $\sum_{k=1}^\infty\alpha_k=\infty$ and 
$\sum_{k=1}^\infty \alpha_k^2 <
  \infty.$ Then, if the optimal set $X^*$
  is nonempty, the iterate
  sequence $\{w_{i,k}\}$ of each agent $i\in V$ converges to the same
  optimal point with probability 1 and in mean square.
\end{theorem}
\begin{proof}
  Observe that the conditions of Lemma~\ref{lemma:biralpha} are
  satisfied. Letting $z = x^*$ for some $x^*\in X^*$, taking conditional
  expectations and using the bounds on the error moments, we obtain
  for any $x^*\in X^*$ and any $k$, with probability 1,
  \begin{align*}
    \sum_{i=1}^{m} \EXP{\|v_{i,k+1} - x^*\|^2 \mid F_k} \leq&
    \sum_{i=1}^{m} \| v_{i,k} - x^*\|^2 - 2 \alpha_{k+1} \l( f(y_k) -
    f^* \r) \\ + &2 \alpha_{k+1} \l(\max_{i \in V} C_i \r)
    \sum_{j=1}^{m} \|y_k - w_{j,k}\|\cr + &2 \alpha_{k+1}
    \sum_{i=1}^{m} \mu_{i,k+1}\|v_{i,k} -x^*\| +\alpha_{k+1}^2
    \sum_{i=1}^{m} \l(C_i + \nu_i \r)^2,
   \end{align*}
  where $f^*=f(x^*),$ and we use the notation 
  $\mu_{i,k+1}=\|\EXP{\epsilon_{i,k+1} \mid F_{k}}\|.$ 
  Using the inequality
  $$2\alpha_{k+1}\mu_{i,k+1}\|v_{i,k} -x^*\| \leq
  \alpha_{k+1}^2\|v_{i,k} -x^*\|^2 + \mu_{i,k+1}^2,$$ 
  we obtain with probability 1, 
  \begin{align}
    \sum_{i=1}^{m} \EXP{\|v_{i,k+1} - x^*\|^2 \mid F_k} \leq&
      \sum_{i=1}^{m} \l( 1 + \alpha_{k+1}^2 \r) \| v_{i,k} - x^*\|^2
      \nn \\&- 2\alpha_{k+1} \l( (f(y_{k}) - f^*) -\l(\max_{i \in V}
      C_i \r) \sum_{j=1}^{m} \|y_k - w_{j,k}\|\r. \nn \\& +
      \sum_{i=1}^{m} \mu^2_{i,k+1} - \l. \frac{1}{2}\,\alpha_{k+1}
      \sum_{i=1}^{m} \l( C_i + \nu_i \r)^2\r). \label{eqn:d4}
   \end{align}
  By Theorem~\ref{thm:asconsensus}, we have with probability 1,
  \[
  \sum_{k} \alpha_{k+1} \| w_{j,k} - y_{k}\| < \infty.
  \]
  Further, since $\sum_{k} \mu^2_{i,k} < \infty$ and $\sum_{k}
  \alpha_{k}^2 < \infty$ with probability 1, the relation in
  (\ref{eqn:d4}) satisfies the conditions of
  Theorem~\ref{lemma:asm}. We therefore have
  \begin{equation}
    \sum_{k} \alpha_{k} (f(y_{k}) - f^*) < \infty, \label{eqn:fin}
  \end{equation}
  and $\|v_{i,k} - x^*\|$ converges with probability 1 and in mean
  square. In addition, by Theorem~\ref{thm:asconsensus}, we have
  $\lim_{k\to\infty}\|w_{i,k} - y_{k}\| =0$ for all $i$, with
  probability 1.  Hence, $\lim_{k\to\infty} \|v_{i,k} -y_{k}\| \to 0$
  for all $i$, with probability 1.  Therefore, $\|y_{k} - x^*\|$
  converges with probability 1 for any $x^*\in X^*.$ Moreover, from
  (\ref{eqn:fin}) and the fact that $\sum_{k} \alpha_{k} = \infty,$ by
  continuity of $f$, it follows that $y_{k},$ and hence $w_{i,k},$
  must converge to a vector in $X^*$ with probability 1 and in mean
  square.
\end{proof}

Note that the result of Theorem~\ref{thm:dss} holds without assuming 
compactness of the constraint set $X$. 
This was possible due to the assumption that both 
the stepsize $\a_k$ and the norms 
$\|\EXP{\epsilon_{i,k+1} \mid F_{k}}\|$ of the conditional errors
are square summable.
In addition, note that the result of Theorem~\ref{thm:dss} remains
valid when the condition
$\sum_{k=0}^\infty\|\EXP{\epsilon_{i,k+1} \mid F_{k}}\|^2<\infty$ for all $i$
is replaced with 
$\sum_{k=0}^\infty\a_{k+1}\|\EXP{\epsilon_{i,k+1} \mid F_{k}}\|<\infty$
for all $i$.

\section{Implications}\label{sec:implications}
The primary source of stochastic errors in the subgradient evaluation
is when the objective function is not completely known and has some
randomness in it.  Such settings arise in sensor network applications
that involve distributed and recursive estimation \cite{Sundhar07b}.

Let the function $f_i(x)$ be given by 
$f_i(x)= \EXP{g_i(x,R_i)},$ where $R_i$ is a random
variable whose statistics are independent of $x.$ The statistics of
$R_i$ are not available to agent $i$ and hence the function $f_i$ is
not known to agent $i.$ Instead, agent $i$ observes samples of $R_i$ in time.
Thus, in a subgradient algorithm for minimizing the 
function, the subgradient must be suitably approximated using the
observed samples.  In the Robbins-Monro stochastic approximation
\cite{Polyak87}, the subgradient $\nabla f_i(x)$ is approximated by
$\nabla g_i(x,r_{i}),$ where $r_{i}$ denotes a sample of $R_i.$ The
associated distributed Robbins-Monro stochastic optimization algorithm
is
\begin{align}
  w_{i,k+1} = P_{X} \l[ v_{i,k} - \alpha_{k+1} \nabla
    g_i\l(v_{i,k},r_{i,k+1}\r) \r],
  \label{eqn:sagrad1}
\end{align}
where $r_{i,k+1}$ is a sample of $R_i$ obtained at time $k.$ 
The expression for the error is
$$ \epsilon_{i,k+1} = \nabla g_i(v_{i,k},r_{i,k+1}) - \EXP{\nabla
g_i(v_{i,k},R_i)}.
$$ If the samples obtained across iterations are independent then
$$\EXP{\epsilon_{i,k+1} \mid F_{k}} = \EXP{\epsilon_{i,k+1} \mid
  v_{i,k}} = 0.$$ If in addition, $\Var{ \nabla g_i(x,R_{i})}$ is
  bounded for all $x \in X$ then the conditions of
  Theorems~\ref{thm:constantstep},~\ref{thm:averi} and \ref{thm:dss}
  are satisfied.

Let us next consider the case when $f_i(x) = \EXP{g_i(x,R_i(x))},$
where $R_i(x)$ is a random variable that is parameterized by $x.$ To
keep the discussion simple, let us assume that $x \in \Re.$ As in
the preceding case, the statistics of $R_i(x)$ are not known to agent
$i,$ but the agent can obtain samples of $R_i(x)$ for any value of $x.$ 
In the
Kiefer-Wolfowitz approximation \cite{Polyak87},
$$ \nabla f_i(x) \approx \frac{g_i\l(x,r_i(x+\beta)\r) - g_i\l(x,r_i(x)
\r) }{\beta},
$$ where $r_i(x)$ is a sample of the random variable $R_i(x).$
The corresponding distributed optimization algorithm is
\begin{align*}
  w_{i,k+1} = P_{X} \l[ v_{i,k} - \alpha_{k+1}
    \frac{g_i\l(v_{i,k},r_i(v_{i,k}+\beta_{i,k+1})\r) -
      g_i\l(v_{i,k},r_i(v_{i,k}) \r) }{\beta_{i,k+1}}
    \r],
\end{align*}
where $\beta_{i,k+1}$ is a positive scalar. In this case, the error is
\[
\epsilon_{i,k+1} = \frac{g_i\l(v_{i,k},r_i(v_{i,k} +\beta_{i,k+1}) \r) -
    g_i\l(v_{i,k},r_i(v_{i,k})\r) }{\beta_{i,k+1}} - \nabla
    f_i(v_{i,k}).
\]
If the function $g_i$ is differentiable then $\EXP{\epsilon_{i,k+1} \mid
v_{i,k}}$ is of the order $\beta_{i,k+1}.$ Thus, the conditions
on the mean value of the errors can be controlled through the sequence
$\{\beta_{i,k}\}$ and the conditions in
Theorems~\ref{thm:constantstep},~\ref{thm:averi} and \ref{thm:dss} can
be met by suitably choosing the sequence $\{\beta_{i,k}\}$.

\section{Discussion}
\label{sec:conclusion}
We studied the effects of stochastic subgradient errors on
distributed algorithm for network of agents with time-varying connectivity. 
We first considered very general errors
with bounded second moments and obtained explicit bounds on the
agent disagreements and on the
expected deviation of the limiting function value from the optimal. The bounds
are explicitly given 
as a function of the network properties, objective function and
the error moments. For networks that are connected at all times and $\eta$ is independent of the size of the network, the
bound scales as $\alpha \left(\max_{i\in V}\{C_i+ \nu_i\}\right)^2 m^4,$ 
where $m$ is the number of
agents in the network, $\alpha$ is the stepsize limit, and 
$C_i$ and $\nu_i^2$ are respectively 
the subgradient norm bound and the bound on the second 
moment of the subgradient errors for agent $i$.  For the constant
stepsize case, we obtained a bound on the performance of the algorithm
after a finite number of iterations. There, we showed that
deviation from the ``error-bound'' diminishes at rate 
$\frac{1}{t},$ where $t$ is the number of iterations.   
Finally, we proved that when
the expected error and the stepsize converge to $0$ sufficiently fast, 
the agents reach a consensus and the iterate sequences of agents 
converge to a common optimal point with probability 1 and in mean square.

We make the following remarks. First, it can be shown that the 
disagreement results
in Corollary~\ref{cor:expconsensus} and Theorem~\ref{thm:asconsensus}
hold even when the agents use non-identical stepsizes. However, with
non-identical agent stepsizes there is no guarantee that the sum of
the objectives rather than a weighted sum, is minimized. 

Future work includes several important extensions of the distributed
model studied here.  At first, we have assumed no communication delays
between the agents and synchronous processing. An important extension
is to consider the properties of the algorithm in asynchronous
networks with communication delays, as in \cite{Tsitsiklis84}.  At
second, we assumed perfect communication scenario, i.e., noiseless
communication links. In wireless network applications, the links are
typically noisy and this has to be taken into consideration.  At
third, we have considered the class of convex functions. This
restricts the number of possible applications for the algorithm.
Further research is to develop distributed algorithms when the
functions $f_i$ are not convex.

\bibliographystyle{siam} 
\bibliography{dssm}
\end{document}